\documentclass{article}
\usepackage{amssymb, amsmath,epsfig}
\usepackage{color}

\newcommand{\mc}[1]{{}}  

\newcommand \na {a}
\newcommand \nb {b}

\newcommand{\hol}{\operatorname{hol}}

\newcommand {\mm} {{\bf m}}

\newcommand \K {{\bold C}}

\DeclareMathRadical{\sqrtsign}{symbols}{"70}{largesymbols}{"70}
\newcommand{\bb}{\mathbb}
\newcommand{\gothic}{\mathfrak}

\newcommand{\half}{{\bb H}}

\newcommand{\natls}{{\bb N}}

\newcommand{\reals}{{\bb R}}

\newlength{\figboxwidth}             
\newcommand{\makefig}[3]{
        \begin{figure}[htb]
        \refstepcounter{figure}
        \label{#2}
        \begin{center}
                #3~\\
                \smallskip
                Figure \thefigure.  #1
        \end{center}
        \medskip
        \end{figure}
}

\renewcommand{\bold}[1]{\medskip \noindent {\bf #1 }\nopagebreak}

\newcommand{\qed}{\hfill $\Box$}
\newcommand{\isom}{\cong}

\newcommand{\cross}{\times}
\newcommand{\st}{\;\: : \;\:}         %Such that

\newcommand{\Vol}{\operatorname{Vol}}

\def\@ifundefined#1#2#3%
  {\expandafter\ifx\csname#1\endcsname\relax#2\else#3\fi}

\@ifundefined{theoremstyle}{
}{
\theoremstyle{plain} %default
}
\newtheorem{theorem}{Theorem}[section]
\newtheorem{prop}[theorem]{Proposition}
\newtheorem{proposition}[theorem]{Proposition}
\newtheorem{lemma}[theorem]{Lemma}

\newtheorem{corollary}[theorem]{Corollary}

\@ifundefined{theoremstyle}{
}{
\theoremstyle{definition} %default
}

\newcommand{\cC}{{\cal C}}

\newcommand{\cF}{{\cal F}}
\newcommand{\cG}{{\cal G}}

\newcommand{\cK}{K}
\newcommand{\cKM}{{\rm K}}
\newcommand{\cKT}{\widetilde{\rm K}}
\newcommand{\cKQ}{{\cal K}}

\newcommand{\cM}{{\cal M}}
\newcommand{\cN}{{\cal N}}

\newcommand{\cP}{{\cal P}}
\newcommand{\cQ}{{\cal Q}^{1}}

\newcommand{\cT}{{\cal T}}

\newcommand{\cX}{{\cal X}}

\newcommand{\cl}{\ell_{min}}
\newcommand{\ee}{{\rm E}}
\mathchardef\GG="321D
\renewcommand{\gg}{{\gothic g}}

\newcommand{\gp}{{\gothic p}}

\newcommand{\te}{\mathcal{T}}
\newcommand{\vo}{\operatorname{Vol}}

\newcommand{\RR}{\operatorname{Re}}   
\newcommand{\II}{\operatorname{Im}}

\newcommand{\tw}{\operatorname{tw}}

\newcommand{\ch}{\S}

\newcommand{\e}{{\operatorname{Ext}}}

\newcommand{\Nbhd}{\operatorname{Nbhd}}
\title{Counting closed geodesics in Moduli space}
\author{Alex Eskin\thanks{partially supported by NSF grants
    DMS 0244542 and DMS 0604251 and the Clay foundation} 
and Maryam Mirzakhani\thanks{partially supported by the Clay foundation }}
\begin{document}

%1)  \gamma was too many things- called the mapping class group elements \gg-
% gamma is reserved for the parametrized geodesics--- and \gamma(L)
% 2) uniform using \cQ\cM====> adding 1 in the case of unit area
%3)  Mod and \gamma---------------\cN and \Delta uniform notation
% 4) all capital X and Y
% 5) the map p is called \gp to be easier to keeop track of
%6) Added a part on extremal lengths/Minsky;s product theorem, more notation....
%7) Moved Lemma 2.3 in the background section and added a short proof for Lemma2
%8) added a ref. for Maskit's bound on extremal length
%9) changing the notation for \asymp===>\approx for uniform use
% uniform use of d_{\cT} for Teich distance
% the length of the shortest saddle connection:::\cl the referee didnt like \ell----
%The same applies to
%the statement that up to uniformly bounded index, the set IX,L only consists of
%twists about the curves ?????????????????
%OCT 12:More typos-

\maketitle
\section{Introduction}
%MAYBE P_J} SHOULD BE P_{J,\TAU}?????????????OTHERWISE LEMMA 4.6 WOULD BE REALLY WEIRED
%BIG CHANGE IN LEMMA 4.6/ NOT USING G HERE
Let $\cM_g$ denote the moduli space of closed Riemann surfaces of genus $g$.
We may write $\cM_g = \cT_g/\Gamma_{g}$, where $\cT_g$ is the
Teichm\"uller space of genus $g$ surfaces, and $\Gamma_{g}$ is the mapping
class group.
In this paper, we investigate properties of Teichm\"uller geodesics on $\cM_{g}.$ 
 Let $N(R)$ denote the number of closed Teichm\"uller geodesics in
$\cM_g$ of length at most $R$. Then $N(R)$ is also the number of 
conjugacy classes of pseudo-Anosov elements of the mapping class group
of translation length at most $R$. Our main result is the following:
\begin{theorem}
\label{theorem:asymp:allgeodesics}
As $R \to \infty$, we have
\begin{displaymath}
N(R) \sim \frac{e^{hR}}{hR},
\end{displaymath}
where $h = 6g-6$.
\end{theorem}
In the above theorem and below, the notation
$A \sim B$ means that the ratio $A/B$ tends to $1$. 
Even though we assume that the surface has no punctures here,
most of the results also hold on $\cM_{g,n}$, provided $h$ is replaced
by $6g-6+2n$. 

In the proof of Theorem~\ref{theorem:asymp:allgeodesics}, the key is
to estimate the number of closed geodesics which stay outside of compact
sets. 
Let $N_j(\delta,R)$ denote the number of closed geodesics $\gamma$ of
length at most $R$ in $\cM_g$ such that for each point $X \in \gamma$,
$X$ has at least $j$ simple closed curves of hyperbolic length less than $\delta$. 

\begin{theorem}
\label{theorem:having:j:short:curves}
Given  $\epsilon > 0$ there exist $\delta > 0$ and $C = C(\epsilon)$
such that for all $j$, $1\le j \le 3g-3$, and all $R>0$,
\begin{displaymath}
N_j(\delta,R) \le C e^{(h-j+\epsilon)R}.
\end{displaymath}
\end{theorem}
Let $\cKM$ be a compact subset of $\cM_g$. We let $N^{\cKM}(R)$ denote
the number of geodesics in $\cM_g$ of length less than $R$ which never
intersect $\cKM$. Letting $j=1$ in
Theorem~\ref{theorem:having:j:short:curves}, we obtain the following:
\begin{corollary}
\label{cor:number:staying:outside:compact}
For every $\epsilon > 0$ there exist a compact set $\cKM \subset \cM_g$ and $C =
C(\epsilon)>0$ such that for all $R > 0$, 
\begin{displaymath}
N^{\cKM}(R) \le C e^{(h-1+\epsilon)R}.
\end{displaymath}
\end{corollary}
Corollary~\ref{cor:number:staying:outside:compact} is complimentary to
the following result:
\begin{theorem}[Rafi, Hammenstadt]
\label{theorem:rafi:hammenstadt}
For any compact $\cKM \subset \cM_g$, and sufficiently large $R$,
\begin{displaymath}
N^{\cKM}(R) \ge e^{(h-1)R}.
\end{displaymath}
\end{theorem}

\begin{subsection}{Previous Results}
The first results on this problem are due to Veech \cite{Veech}. He
proved that there exists a constant $c_2$ such that
\begin{displaymath}
h \le \liminf_{R \to \infty} \frac{\log N(R)}{R} \le \limsup_{R
  \to \infty} \frac{\log N(R)}{R} \le c_2
\end{displaymath}
and conjectured that $c_2 = h$. In a remark in a paper by Ursula
Hamenstadt \cite{Hamenstadt:Bernoulli} (see also
\cite{Hamenstadt:dynamics}), in
which the main focus is different, she proves that $c_2 \le (6g - 6 +
2n)(6g - 5 + 2n)$. 

Sasha Bufetov \cite{Bufetov} proved the formula
\begin{equation}
\label{eq:cover:asymp}
\lim_{R \to \infty} \frac{\log \tilde{N}(R)}{R} = h, 
\end{equation}
where $\tilde{N}(R)$ is the number of periodic orbits of the
Rauzy-Veech induction such that the log of the norm of the
renormalization matrix is at most $R$. This is a closely related
problem; essentially $\tilde{N}(R)$ counts closed geodesics on a
certain finite cover of $\cM_g$. However the equation
(\ref{eq:cover:asymp}) does not easily imply
\begin{equation}
\label{eq:log:asymp}
\lim_{R \to \infty} \frac{\log N(R)}{R} = h. 
\end{equation}
Very recently, Kasra Rafi \cite{Rafi} proved
Corollary~\ref{cor:number:staying:outside:compact} (which implies
(\ref{eq:log:asymp})) for the case of the five-punctured sphere.

 We note that (\ref{eq:log:asymp}) is an immediate consequence of
Theorem~\ref{theorem:asymp:allgeodesics}, which is a bit more precise.
In order to prove Theorem~\ref{theorem:asymp:allgeodesics} one needs
Corollary~\ref{cor:number:staying:outside:compact} and certain
recurrence results for geodesics, which are based on \cite{Jayadev}.

\bold{Remarks.}
\begin{itemize}
\item The problem of understanding the asymptotics 
of the number $N_{M}(R)$ of primitive closed geodesics of length $\leq R$ 
on a given manifold $M$ has been investigated intensively. An asymptotic formula for $N_{M}(R)$ on a compact hyperbolic surfaces was first proved by Huber. See $\ch 9$ in \cite{B:book} and references within for related work of Hejhal, Randol and Sarnak.  

More generally, Margulis proved that on a compact $n$-manifold $M$ of negative curvature 
$$N_{M} (R)\sim \frac{e^{hR}}{h R},$$
where $h$ is the topological entropy of the geodesic flow. 
In this case the techniques from uniformly hyperbolic dynamics can be applied to study the geodesic flow on the unit tangent bundle of $M$.
See \cite{Margulis:thesis} for the proof of Margulis' theorem, and related results on Anosov and hyperbolic flows.  

\item The main difficulty for proving Theorem \ref{theorem:asymp:allgeodesics} is the fact that the Teichm\"uller flow is not hyperbolic. 
In order to overcome this difficulty, first in \ch \ref{sec:thin:part} we use Minsky's product region theorem \cite{Minsky} to prove that the geodesic flow is biased toward the 
{\it thick part} of the moduli space. This result implies Theorem \ref{theorem:having:j:short:curves}.
Then we use the basic properties of the Hodge norm \cite{ABEM} to prove a closing lemma for the Teichm\"uller geodesic flow in \ch \ref{closing}.
But the Hodge norm behaves badly near smaller strata, i.e. near points with degenerating zeros of the quadratic differential. 
However, in \ch \ref{rec} we show that
 the number of closed geodesics $\gamma$ of length at most $R$ such that $\gamma$ spends at least
$\theta$-fraction of the time outside of a compact subset of the principal stratum is exponentially smaller than $N(R)$ (see Theorem \ref{theorem:staying:in:principal:stratum}).
Finally, we obtain the counting result using the fact that in any compact subset of the principal stratum, the geodesic flow is uniformly hyperbolic (\cite{Veech}, \cite{forni}, and \cite{ABEM}). 
\item 
 By results in \cite{Hamenstadt:BM} the normalized geodesic flow invariant measure supported on the set of closed geodesics of length $\leq R$ in $\mathcal{Q}^1\mathcal{M}_{g}$ become equidistributetd with respect to the Lebesgue measure $\mu$ (see \ch \ref{qd}) as $R\rightarrow \infty$. 

\item In a forthcoming joint work with Kasra Rafi, we generalize the results in this paper to the case of other strata of moduli spaces of Abelian and quadratic differentials.

\end{itemize}
\end{subsection}
\noindent
{\bf Notation.} In this paper, $A \approx B$ means that $A/C< B < A C$ for some universal constant $C$ which only depends genus $g$.
Also, $A=O(B)$ means that $A < B C,$ for some universal constant $C$, which again could depend on $g$. \\

\bold{Acknowledgements.} We would like to thank Kasra Rafi who
suggested a major simplification of the proof of
Theorem~\ref{theorem:having:j:short:curves}. We would also like to
thank Howard Masur for many useful discussions (including the
proof of Lemma~\ref{lemma:uniform:vol:ball}) and also
William Cavendish whose lecture notes were used for part of
this paper. We would like to thank the referee for helpful comments.
\section{Background and notation}
\label{sec:back}
In this section, we recall definitions and known results about the Teichm\"uller geometry 
of $\mathcal{M}_{g}.$ For more details see \cite{Hubbard}.

\subsection{Teichm\"uller space} 
A point in the {\it Teichm\"uller space} $\cT_{g}$
is a complex curve $X$ of genus $g$ 
equipped with a diffeomorphism $f: S_{g} \rightarrow X$. The map $f$ provides a {\it marking} on $X$ by $S_{g}$. Two
marked surfaces $f_{1}:\; S_{g} \rightarrow X$ and $f_{2}:\; S_{g} \rightarrow Y$
define the same point in $\cT_{g}$ if and only if $f_{1} \circ f_{2}^{-1} : Y \rightarrow X$
is isotopic to a holomorphic map.
By the uniformization theorem, each point $X$ in $\cT_{g}$ has a metric of constant 
curvature $-1$. 
The space $\cT_{g}$ is a complex manifold of dimension $3g-3$, diffeomorphic to a cell. 
Let $\Gamma_{g}$
denote  the mapping class group of $S_{g}$, or in other words
the group of isotopy classes of orientation preserving self
homeomorphisms of $S_{g}$. The mapping
class group $\Gamma_{g}$ acts
on $\cT_{g}$ by changing the marking. The quotient space
$$\cM_{g}= \cT_{g}/\Gamma_{g}$$
is the moduli space of Riemann surfaces homeomorphic to $S_{g}$.
The space $\cM_{g}$ is an orbifold. However, it is finitely covered by a manifold, and  
$$\pi_{1}^{orb}(\cM_{g})=\Gamma_g.$$
Let $$\gp: \cT_g \rightarrow \cM_g$$ denote the natural map from $\cT_g$ to $\cM_g= \cT_g/\Gamma_{g}$. 
The
Teichm\"uller metric on marked surfaces is defined by
 $$d_{\cT}((f_{1} : S_{g} \rightarrow X_{1}), (f_{2}: S_{g} \rightarrow X_{2})) = \frac{1}{2} \inf  (\log K(h)),$$ 
where $h : X_{1} \rightarrow X_{2}$ ranges over all quasiconformal maps isotopic to $f_{1}o f_{2}^{-1}.$ Here $K(h) \geq 1$ is the dilatation of $h$.\\
\noindent
{\bf Dilatation of pseudo-Anosov elements.} According to the Nielsen-Thurston classification, every irreducible mapping class element $\gg \in \Gamma_{g} $ of
infinite order has a representative which is a pseudo-Anosov homeomorphism \cite{Th}. 
By a theorem of Bers, every closed geodesic in $\mathcal{M}_{g}$ is the unique loop of minimal length in its homotopy class.
Given a pseudo-Anosov $\gg \in \Gamma_{g}$ the dilatation of $\gg$ is defined by $K(\gg)$.
Then $\log(K(\gg))$ is the translation length of $\gg$ as an isometry of $\cT_{g}$ \cite{Bers}. In other words, 
$$\mathcal{L}(S_{g})=\{ \log(K(\gg))\; | \gg \in \Gamma_{g} \; \mbox {pseudo-Anosov}\;\} $$
is the length spectrum of $\cM_{g}$ equipped with the Teichm\"uller metric.\\
By \cite{AY} and \cite{Iv} $\mathcal{L}_{g}$ is discrete subset of ${\Bbb R}$. 
Hence the number $N(R)$ of conjugacy classes of pseudo-Anosov elements of the group $\Gamma_{g}$ with dilatation factor $K(\gg) \leq R$ is finite.
We remark that for any pseudo-Anosov $\gg \in \Gamma_{g}$ the number $K(\gg)$ is an algebraic number. Moreover $\log(K(\gg))$ is equal to the minimal topological entropy of any element in the same homotopy class \cite{FLP}. \\  

\subsection{Moduli space of quadratic differentials}\label{qd}
The cotangent space of $\cT_{g}$ at a point $X$ can be identified with the vector space $Q(X)$ of holomorphic quadratic differentials on $X.$ 
Recall that given $X \in \cT_{g}$, a quadratic differential $q \in Q(X)$ is a tensor locally given by $\phi(z) dz^{2}$ where $\phi$ is holomolphic. 
Then the space $\mathcal{Q}\te_{g}= \{(q,X) \;|\; X \in \te_{g}, q \in Q(X)\}$ is the cotangent space of 
$\te_{g}.$
In this setting, the Teichm\"uller metric corresponds to the norm 
$$\parallel q \parallel_{\cT} = \int _{X} |\phi(z)|\; |dz|^{2}$$
on $\mathcal{Q}\te_{g}.$
Let $\mathcal{Q}\cM_{g}\cong \mathcal{Q}\te_{g}/ \Gamma_{g}$. 
Finally, let $\cQ\cT_g$ denote the Teichm\"uller 
space of unit area quadratic differentials on surfaces of genus $g$, and 
$\mathcal{Q}^{1}\cM_{g}\cong \cQ\te_{g}/ \Gamma_{g}.$
For simplicity, we let $\pi$ denote 
the natural projection maps 
 $$\pi: \cQ\cM_g \to \cM_g,$$
 and  
 $$\pi: \cQ\cT_g \to \cT_g.$$
Although the value of $q\in \cQ(X)$ at a point $x \in X$ depends on the local coordinates, the zero set of $q$ is well defined.
As a result, there is a natural stratification of the space $\mathcal{Q}\mathcal{M}_{g}$ by the multiplicities of zeros of $q$.
Define $\mathcal{Q}\mathcal{M}_{g}( a_{1},\ldots,a_{k}) \subset \mathcal{Q}\mathcal{M}_{g}$ to be the subset consisting of pairs $(X,q)$ of holomorphic quadratic differentials on $X$ with $k$ zeros with multiplicities $(a_{1},\ldots,a_{k})$. 
Then $$\mathcal{Q}\mathcal{M}_{g}= \bigsqcup\limits_{(a_{1},\ldots,a_{k})} \mathcal{Q}\mathcal{M}_{g}(a_{1},\ldots,a_{k}).$$
It is known that each $ \mathcal{Q}\mathcal{M}_{g}(a_{1},\ldots,a_{k})$ is an orbifold of dimension $4g-4+2k$.
In particular $ \dim(\mathcal{Q}\mathcal{M}_{g}(1,\ldots,1))=\dim(\mathcal{Q}\mathcal{M}_{g}).$

We recall that when $g>1,$ the Teichm\"uller metric is not even Riemannian. However, geodesics in this metric are well understood. 
A quadratic differential $q \in \mathcal{Q}\te_{g}$ with zeros at $x_{1},\ldots x_{k}$ is determined by an atlas of charts $\{ \phi_{i}\}$ mapping open subsets of $S_{g}-\{x_{1},\ldots,x_{k}\}$ to ${\Bbb R}^{2}$ such that the change of coordinates are of the form $v \rightarrow \pm{v}+c.$
Therefore the group $\operatorname{SL}_{2}({\Bbb R})$ acts naturally on $\mathcal{Q}\cM_{g}$
by acting on the corresponding atlas; given $A \in \operatorname{SL}_{2}({\Bbb R})$, $A \cdot q \in \mathcal{Q}\cM_{g}$ is determined by the new atlas $\{ A \phi_{i}\}.$
The action of the diagonal subgroup
$g_{t}=\begin{bmatrix} e^{t/2}&0\\
0 & e^{-t/2} \end{bmatrix}$
is the Teichm\"uller geodesic flow for the Teichm\"uller metric.
Moreover we have \cite{Veech},\cite{Masur1}:
\begin{theorem}{({\bf Veech-Masur})}\label{VM}
The space $\cQ\cM_{g} $ carries a unique probability measure $\mu=\mu_{g}$ in the Lebesgue measure class such that :
\begin{itemize}
\item the action of $SL_{2}({\Bbb R})$ is volume preserving and ergodic;
\item Teichm\"uller geodesic flow is mixing.
\end{itemize}
\end{theorem}
By this theorem, $\cQ\cT_g$ carries a natural normalized smooth measure $\mu=\mu_{g}$, 
preserved by the action of $\Gamma_{g}$. We set 
\begin{equation}\label{measure}
\mm=\pi_*\mu.
\end{equation}
\noindent
{\bf Remark.}
In fact, the Teichm\"uller flow on $\cQ\cM_{g}$ is exponentially mixing with respect 
to $\mu$ \cite{Avila:Resende}. See also \cite{AGY:mixing}. However, we
will only use Theorem $\ref{VM}$ in this paper.   
\subsection {Period coordinates on the strata}\label{pc}
A {\it saddle connection} on $q \in \cQ\cT_{g}$ is a geodesic segment which joins a pair of singular points without passing
through one in its interior. In general, a geodesic segment $e$ joining two zeros of
a quadratic differential $q=\phi dz^{2}$ determines a complex number $\hol_{q}(e)$ (after choosing a branch of $\phi^{1/2}$ and
an orientation of $e$) by
$$\hol_{q}(e)=\RR(\hol_{q}(e))+\II(\hol_{q}(e)),$$
where
$$\RR(\hol_{q}(e))= \int\limits_{e} \RR(\phi^{1/2}),$$
and
$$ \II(\hol_{q}(e))=\int\limits_{e} \II(\phi^{1/2}).$$
We recall that the period coordinates give $\mathcal{Q}\te_{g}(a_{1},\ldots,a_{k})$ the structure of a piecewise linear manifold.
 For notational simplicity we discuss the case of $\mathcal{Q}\te_{g}(1,\ldots,1).$
Given $q_{0} \in \mathcal{Q}\te_{g}(1,\ldots,1)$ there is a triangulation $E$ of the underlying surface by saddle connections, $h=6g-6$ directed edges $\delta_{1},\ldots,\delta_{h}$ of $E$, and an open neighborhood  
$U_{q_{0}} \subset \mathcal{Q}\te_{g}(1,\ldots,1)$ of $q_{0}$ such that the map 
$$\psi_{E,q_{0}}: \mathcal{Q}\te_{g}(1,\ldots,1) \rightarrow {\Bbb C}^{6g-6}$$
by 
$$\psi_{E,q_{0}}(q)=(\hol_{q}(\delta_{i}))_{i=1}^{h}$$
is a local homeomorphism. Also for any other geodesic triangulation $E^{'}$ the map 
the map $\psi_{E^{'},q_{0}}\circ\psi_{E,q_{0}}^{-1}$ is linear.
For a discussion of these coordinates see \cite{Masur:Smillie}. The measure $\mu_{g}$ in Theorem \ref{VM} (up to a constant) is given by the piecewise linear structure of
$\mathcal{Q}\cM_{g}.$ This measure, up to normalization, also coincides with the measure defined by the Teichm\"uller norm on the unit cotangent bundle of $\cT_{g}.$
This measure is supported on $\mathcal{Q}^{1} \cM_{g}(1,1,\ldots,1))$; that is, $\mu_{g}(\cQ\cM_{g}-\cQ \cM_{g}(1,1,\ldots,1))=0.$\\
%ADD a A REFERENCE)
%GEODESIC FLOW
%EXTREMAL LENGTH RELATION
%ADD HERE\\
%A ONTEICHMULLER DISTANCE\\
%THE LENGTH SPECTRUM IS DISCRETE\\
%G(X) AND SOME PROPERTIES IF X,Y ARE CLOSE THE RATIO OF G(X) AND G(Y) BOUNDED AND INTEGRABLE????
%B(X,R) DEFINE. MOD MAPPING CLASS GROUP\\
%NOTATION ISSUE WITH AYMP...ALEX USES A DIFFERENT NOTATION\\
%P AND PI HOW ARE THEY USED CORRECTLY\\
%DILATATION OF PA REQUIRES SOME EXPLANATION!!\\
\subsection{Extremal and hyperbolic lengths of simple closed curves}\label{boundedpants}
 Given a homotopy class of a simple closed curve $\alpha$ on a topological surface $S_{g}$ and 
 $X \in \cT_{g}$, let $\ell_{\alpha}(X)$ be the length of the unique 
geodesic in the homotopy class of $\alpha$ with respect to the hyperbolic metric on $X$. 
The
{\it extremal length} of a simple closed curve $\alpha$ on $X$ is defined by
\begin{equation}\label{ext}
\e_{\alpha}(X)=\sup_{\rho}
\frac{\ell_{\alpha}(\rho)^{2}}{\operatorname{Area}(X,\rho)},
\end{equation}
where the supremum is taken over all metrics $\rho$ conformally
equivalent to $X$, 
and $\ell_\alpha(\rho)$ denotes the length of $\alpha$ in the metric
$\rho$. \\
Given simple closed curves $\alpha$ and $\beta$ on $S_{g}$, the intersection number $i(\alpha,\beta)$ is the minimum number of points in which representatives of $\alpha$ and $\beta$ must intersect. In general, we have: 
$$i(\alpha, \beta) \leq \sqrt{\e_{\alpha}(X)} \cdot \sqrt{\e_{\beta}(X)}. $$
The following result \cite{Kerckhoff:Asm} relates the ratios of extremal lengths to the Teichm\"uller distance:
\begin{theorem}({\bf Kerckhoff})\label{ker}
Given $X,Y \in \te_{g}$, the Teichm\"uller distance between $X$ and
$Y$ is given by 
$$d_{\cT}(X,Y)=\sup_{\beta}\log\left(\frac{\sqrt{\e_{\beta}(X)}}{\sqrt{\e_{\beta}(Y)}}\right),$$
where $\beta$ ranges over all simple closed curves on $S_{g}.$
\end{theorem} 
\noindent 
 The relationship between the extremal length and hyperbolic length is complicated; in general, by the definition 
 of extremal length
$$ \frac{\ell_{\alpha}(X)^{2}}{ 4 \pi (g-1)} \leq \e_{\alpha}(X).$$
Also for any $X\in \te_{g}$ the extremal length can be extended continuously to the space of measured laminations \cite{Kerckhoff:Asm} such that 
$$\e_{r \cdot \lambda}(X)=r^{2} \e_{\lambda}(X).$$ As a result given $X$ there exists a constant $a_{X}$ such that 
$$\frac{1}{a_{X}} \ell_{\alpha}(X) \leq \sqrt{\e_{\alpha}(X)} \leq a_{X} \ell_{\alpha}(X). $$
However, by \cite{Maskit},
\begin{equation}\label{Maskit}
  \frac{1}{\pi} \leq \frac{\e_{\alpha}(X)} { \ell_{\alpha}(X)} \leq \frac{1}{2} e^{\ell_{\alpha}(X)/2}.
\end{equation}
 Hence as $\ell_{\alpha}(X) \rightarrow 0, $ $$\frac{\ell_{\alpha}(X)}{\e_{\alpha}(X)} \approx1.$$\\
\noindent
{\bf Bounded pants decompositions.}
Recall that by a theorem of Bers, there exists $L_g$ depending only on $g$ such that for any surface $Y \in \cT_g$ there exists a pants decomposition
$\mathcal{P}_{Y}=\{\alpha_{1},\ldots,\alpha_{3g-3}\}$ 
of $Y$ satisfying $\ell_{\alpha_{i}}(Y) \leq L_{g}.$
By $(\ref{Maskit})$ the extremal length of each $\alpha_{i}$ on $Y$ is
bounded from above by $C_g$, where $C_{g}=e^{L_{g}} L_{g}$ is independent of $Y$. 
We call such a pants decomposition a {\it bounded} pants decomposition for $Y.$\\
Fix a small enough number $\epsilon_{0}>0$ such that:\\
\noindent
{\bf 1).} Any two simple closed curves of extremal length $\leq \epsilon_{0}^2$ are disjoint;\\
\noindent
{\bf 2).} any simple closed curve intersecting a simple closed curve of extremal length $\leq \epsilon_0^2$ on a surface of genus $g$ has extremal length $\geq 2 C_{g}.$ \\
 We say $\alpha$ is {\it short} on $X$ if $\e_{\alpha}(X) \leq \epsilon_{0}^{2}.$ Let $\cC_{X}$ be the set of simple closed curves $\beta$ on the surface $X$ with $\e_\beta(Y) \le \epsilon_0^2$. Note that by the definition any bounded pants decomposition of $X$ contains all the elements of $\cC_{X}.$ 

Now, define $G: \cT_{g} \rightarrow {\Bbb R}_{+}$ by
\begin{equation}\label{GG}
G(Y) = 1+ \prod_{\beta \in \cC_{Y}}
\frac{1}{\sqrt{\e_\beta(Y)}}.
\end{equation}
Note that for any bounded pants decomposition $\mathcal{P}_{X}=\{\alpha_{1},\ldots,\alpha_{3g-3}\}$  
$$G(X) \approx \prod_{\alpha \in \mathcal{P}_{X}}
\frac{1}{\sqrt{\e_\alpha(X)}}.$$
Also, if $d_{\cT}(X,Y)=O(1)$ then $G(X) \approx G(Y).$
By the definition, $G$ induces a proper function on $\cM_{g}.$\\
%% IS G INTEGRABLE. Is it important?????\\
\subsection{Estimating extremal length and Minsky's product theorem}\label{extremal}
Let $\alpha_{1},\ldots \alpha_{j}$ be a collection of disjoint simple closed curves on $S_{g},$ and $\epsilon \leq \epsilon_0^2$. Let 
$$P_{\epsilon}(\alpha_{1},\ldots, \alpha_{j})=\{ X \in \cT_{g}\; | \e_{X}(\alpha_{i}) \leq \epsilon\}. $$
Let $\half^2$ denote the upper half plane model of the hyperbolic plane. 
Then using the Fenchel-Nielsen coordinates on $\cT_{g}$, we can define $ \phi_{0}: P_{\epsilon}(\alpha_{1},\ldots, \alpha_{j}) \rightarrow (\half^2)^j $
by
$$\phi_{0}(X)=(\theta_{1}(X), \frac{1}{\ell_{\alpha_{1}}(X)}, \ldots, \theta_{j}(X), \frac{1}{\ell_{\alpha_{j}}(X)}),$$
where $\theta_{i}$ is the twist coordinate around $\alpha_{i}.$
In fact, following Minsky, we get a map
$$\phi: P_{\epsilon}(\alpha_1,\dots,\alpha_j) \rightarrow  (\half^2)^j \cross \cT',$$
where $\cT'$ is the quotient Teichm\"uller
space obtained by collapsing all the $\alpha_i$. The product region theorem \cite{Minsky}
states that for sufficiently small $\epsilon $ the Teichm\"uller metric on $P_{\epsilon}(\alpha_1,\dots,\alpha_j)$
is within an additive constant of the supremum metric on
$(\half^2)^j \cross \cT'$.
More precisely, let 
$d'(\cdot, \cdot)$ denote the supremum metric on $(\half^2)^j \cross
\cT'$. Then :
\begin{theorem}
\label{PRT}
If $\epsilon>0$ is small enough, then there exists a
constant $b > 0$ depending only on the genus such that
for all $X,Y \in P_{\epsilon}(\alpha_1,\dots,\alpha_j)$, 
$$|d_{\cT}(X,Y) - d'(\phi(X),\phi(Y))| < b.$$ 
\end{theorem}

Note that by Theorem $\ref{ker}$ this statement can be rewritten in terms of the ratios of extremal lengths of simple closed curves on
$X$ and $Y$.

 Let $\mathcal{C}_{g}$ denote the set of all multicurves on $S_{g},$ and $\{\alpha_{1},\ldots, \alpha_{3g-3}\}$ be a pants decomposition of $S_{g}.$
Consider the Dehn-Thurston parameterization \cite{Harer:Penner:book}  
$$ DT: \mathcal{C}_{g} \rightarrow ({\mathbb Z}_{+} \times {\mathbb Z})^{3g-3}$$
defined by 
$$ DT(\beta)=(i(\beta,\alpha_{i}), \tw(\beta,\alpha_{i}))_{i=1}^{3g-3},$$
where $i(\cdot, \cdot)$ denotes the geometric intersection number and
$\tw(\beta,\alpha_{i})$ is the twisting parameter of $\beta$ around $\alpha_{i}$. See $\cite{Harer:Penner:book}$ for more details.

The proof of Theorem \ref{PRT} relies on the following estimates for the extremal lengths of arbitrary simple closed curves on a surface,
see \cite{Minsky} (Theorem~5.1, and equation (4.3)):
\begin{theorem}{\bf (Minsky)}
\label{theorem:minsky:ext}
Suppose $Y \in T_g$, and let $\cP = \cP_{Y} =
\{\alpha_{1},\ldots,\alpha_{3g-3}\}$ be a bounded pants decomposition on $Y$.
Then given a simple closed curve $\beta$ on $S_{g}$, $\e_\beta(Y)$ 
is bounded from above and below by  
\begin{equation}\label{es}
\max_{1 \leq j \leq 3g-3}
\left[\frac{i(\beta,\alpha_{j})^{2}}{\e_{\alpha_{j}}(Y)}+
\operatorname{tw}^{2}(\beta,\alpha_{j})\e_{\alpha_{j}}(Y) \right]
\end{equation}
up to a multiplicative constant depending only on $g$.
\end{theorem}
This result gives an upper bound for the Dehn-Thurston coordinates of a simple closed
curve in terms of its extremal length.

% We use Theorem $5.1$ in \cite{Minsky:ext} to get an upper bound for
% $E(y,L)$ in terms of $G(y).$ 
% See also equation $(4.3)$ in \cite{Minsky:ext}.
We remark that the definition of the twist used in equation $(4.3)$ in
\cite{Minsky} is different from the definition we are using
here. We follow the definition used in $\cite{Harer:Penner:book}.$ 
%In fact, one can check that for annulus $A$ embedded in $x$
%$$|t_{A,x}(\beta)-\frac{\tw(\beta,\alpha)}{i(\alpha,\beta)}|< 1.$$
 In terms of our notation, 
\begin{equation}\label{ours}
\tw(h_{\alpha}^{r}(\beta),\alpha)=\tw(\beta, \alpha)+ r \cdot i(\beta,\alpha_{i}),
\end{equation}
where $h_{\alpha} \in \Gamma_{g}$ is the right Dehn twist around $\alpha$.

 Let $\mathcal{P}_{X}=\{\alpha_{1},\ldots,\alpha_{n}\}$ be a bounded pants decomposition of $X \in \cT_{g}$ (see \ch \ref{boundedpants}).
Then:
\begin{itemize}
 \item 
 For any simple closed curve $\beta$ on 
$S_{g}$, 
 \begin{equation}\label{esst}
 \tw(\beta,\alpha_{i})\leq c_{1} \cdot  \frac{\sqrt{\e_{\beta}(X)}}{\sqrt{\e_{\alpha_i}(X)}} , \qquad i(\beta,\alpha_{i})\leq \sqrt{\e_{\beta}(X)}\cdot  \sqrt{\e_{\alpha_i}(X)},
 \end{equation}
where $c_{1}>0$ only depend on $g.$
\item Let $\alpha \in \mathcal{P}_{X}$ and $r_{0}>0.$
We claim that if 
$$d_{\cT} (h_{\alpha}^{k}(X), X) \leq r_{0}, $$
then $$|k| \leq c \; \frac{e^{r_{0}}}{\e_{\alpha}(X)},$$
where $c>0$ only depends on $g$ and
$h_{\alpha} \in \Gamma_{g}$ is the right Dehn twist around $\alpha$. 
Note that by \cite{Minsky}, if $\alpha \in \mathcal{P}_{X}$ then there exists a simple closed curve $\beta$ such that 
\begin{equation}\label{product}
\e_{\beta}(X) \approx \frac{1}{\e_{\alpha}(X)},
\end{equation}
and $i(\alpha,\beta)=1.$
We can use Theorem $\ref{theorem:minsky:ext}$ to estimate $\sqrt{\e_{h_{\alpha}^{-k}(\beta)}(X)}= \sqrt{\e_{\beta}(h_{\alpha}^{k}(X))}.$ 
As a result 
$$ \frac{1}{\sqrt{\e_{\alpha}(X)}}+ |k| \sqrt{\e_{\alpha}(X)} \leq  c' \; \sqrt{\e_{h_{\alpha}^{-k}(\beta)}(X)},$$ 
where $c'$ only depends on $g.$
By Theorem $\ref{ker}$, we have: 
$$ \sqrt{\e_{h_{\alpha}^{-k}(\beta)}(X)} \leq e^{r_{0}} \sqrt{\e_{\beta}(X)} .$$
Now from (\ref{product}) we get 
\begin{equation}\label{simplebound} 
\frac{1}{\sqrt{\e_{\alpha}(X)}}+ |k| \sqrt{\e_{\alpha}(X)}= O(e^{r_{0}} \sqrt{\e_{\beta}(X)}) \Longrightarrow |k|=O\left( \frac{e^{r_{0}}}{\e_{\alpha}(X)}\right). 
\end{equation}
In other words, the number of twists around $\alpha_i \in \mathcal{P}_{X}$ which one
can take and still stay in $B_{\cT}(X,r)$ is $O(e^{r}/a_i), $ where $a_{i}=\e_{\alpha_{i}}(X).$
A generalization of this argument is used in the proof of Lemma \ref{Mmain} \cite{ABEM}.\\
\noindent
{\bf Remark.}
As a result of Theorem $\ref{theorem:uniform:upper:bound}$, for any $X \in \cT_{g}$ the number of $\gg \in \Gamma_{g}$ such that $d_{\cT} (\gg \cdot X, X) \leq r$ is 
$O( \prod_{i=1}^{3g-3} r\;e^{r}/a_i).$ 
When $r>0$ is fixed, one can obtain this bound by showing that up to uniformly bounded index, the set of $\{\gg\in \Gamma_{g} \; |\; d_{\cT}(\gg\cdot X,X) \leq r\} $ only consists of
twists about $\alpha_{i} \in \mathcal{P}_{X}$. 
 \end{itemize}
\section{ Net points in Teichm\"uller space} 
Recall that a set $\cN$ is a $(c_{1},c_{2})$ separated net on a metric space $\cX$
if $\cN \subset \cX$, every point of $\cX$ is within $c_2$ of
a net point, and the minimal distance between net points is at least
$c_1$. 
In \ch \ref{sec:thin:part} we will use nets in order to discretize $\cT_{g}.$

 \subsection{Volumes of balls of fixed radius} 
Let $B_{\cT}(X,L) \subset \cT_{g}$ denote the ball of radius $L$ with respect to the Teichm\"uller metric, and let $\mm$ be the smooth measure defined by $(\ref{measure})$ on $\cT_{g}$. Then we have:
\begin{lemma}
\label{lemma:uniform:vol:ball}
There exists $L_0 > 0$ (depending only on $g$) and 
for every $L > L_0$ there exist constants $0 < c_1 < c_2$ 
such that for all $X \in \cT_g$,
\begin{displaymath}
c_1 \le \mm(B_{\cT}(X,L)) \le c_2
\end{displaymath}
The constants $c_1$ and $c_2$ depend on $L$ and $g$ (but not on $X$). 
\end{lemma}

\bold{Proof.} Suppose $X \in \cT_g$. Let $\alpha_1, \dots, \alpha_k$
be the simple closed curves on $X$ with extremal length less than $C_{g}^2$.
Let $a_i$ denote the extremal length of $\alpha_i$. Note that by $(\ref{Maskit})$ the hyperbolic length
$\ell_{\alpha_{i}}(X)\approx a_i$. 
By Theorem $\ref{ker}$ for all $Y \in B_{\cT}(X,L)$, 
\begin{displaymath}
\e_Y(\alpha_i) \le L_{1}^2\; \e_X(\alpha_i) \le C\; L_{1}^2 a_i, 
\end{displaymath}
where $L_{1}=e^{L}$, and $C$ depends only on the genus. Then,
in view of the definition of extremal length (equation (\ref{ext})), for any area $1$
holomorphic quadratic differential $q$ on $Y \in B_{\cT}(X,L)$, 
\begin{equation}
\label{eq:flat:extremal:hyperbolic}
\ell_q(\alpha_i)^2 \le \e_Y(\alpha_i) \le C L_{1}^2 a_i,
\end{equation}
where $\ell_q(\cdot)$ denotes length in the flat metric defined
by $q$, and $C$ depends only on the genus. Thus any flat metric
in the conformal class of a surface in $B_{\cT}(X,L)$ has closed curves of
flat length at most $L_{1} \sqrt{C a_i}$. 

Let $\cF$ be a fundamental domain for the action of $\Gamma_{g}$ on
$\cT_g$. For $k \le 3g-3$ and $0 \le c_1 \le c_2 \le \dots c_k \le 1$,
let $Q(c_1, \dots, c_k)$ denote the subset of $X \in \cT_g$ for which
there exist disjoint curves $\beta_1, \dots, \beta_k$ with
$\e_X(\beta_i) \le c_i$. Then, for any area $1$ holomorphic
quadratic differential $q \in \pi^{-1}(X)$, the length of $\alpha_i$
in the flat metric induced by $q$ is at most $\sqrt{c_i}$. 
Then, by the definition of the measure $\mm(\cdot)$, 
\begin{equation}
\label{eq:measure:short:curves}
\mm(Q(c_1, \dots, c_k) \cap \cF) \le C \prod_{i=1}^k c_i.
\end{equation}
In view of (\ref{eq:flat:extremal:hyperbolic}), 
\begin{displaymath}
B_\cT(X,L) \subset Q(C L_1^2 a_1, \dots, C L_1^2 a_k). 
\end{displaymath}
Then, in view of (\ref{eq:measure:short:curves}), for any $\gg \in
\Gamma_{g}$, 
\begin{equation}
\label{eq:measure:in:fund:domain}
\mm(B_{\cT}(X,L)\cap \gg \cF) \le C \prod_{i=1}^k C L_{1}^2 a_i.
\end{equation}
Let $I_{X,L}$ denote the set of elements $\gg \in \Gamma_{g}$ such 
that $B_{\cT}(X,L) \cap \gg \cF$ is non-empty. We will now estimate
the size of $I_{X,L}$. 
%PLEASE EXPLAIN HERE!!!!!!!!!!!!!!!!!!!!!!!!!!!!!!!!!!!!!!!!!!!!!!!!!!!!!!!!!!!!!!!!!!!!!!!!!!!!!!!!!

%One can check that up to uniformly bounded index (the bound depending on $g$
%and $L$) $I_{X,L}$ consists only of twists around the $\alpha_i$. 

Note that as $\cM_{g}$ is finitely covered by a manifold, (\ref{simplebound}) implies that 
$$
|I_{X,L}| \leq C_{L} G(X)^{2},$$
where $$G(X)\approx \prod_{i=1}^k \frac{1}{\sqrt{a_i}},$$ and $C_L$ only depends on $L$ and $g.$
Also, See Theorem \ref{theorem:uniform:upper:bound} (for the case of $\mathcal{B}=\emptyset$).
 
 On the other hand, by Theorem \ref{theorem:minsky:ext}, the number of twists around $\alpha_i$ which one
can take and still stay in $B_{\cT}(X,L)$ is $\approx L_{1}^2/a_i$. Thus, 
\begin{equation}
\label{eq:IXL}
|I_{X,L}| \approx C'_L \prod_{i=1}^k \frac{1}{a_i}. 
\end{equation}
See \ch \ref{extremal}. The upper bound of Lemma~\ref{lemma:uniform:vol:ball} follows
from (\ref{eq:measure:in:fund:domain}). 

We now briefly outline the proof of the lower bound. Let $R \subset \cQ\cM_g/\Gamma_{g}$ be the
set of flat structures such that each $S \in R$ has flat cylinders
$C_i$ with width (i.e. core curve) between $\sqrt{a_i}$ and
$\sqrt{a_i}/2$ 
and height between $1/(g \sqrt{a_i})$ and $1/(2g \sqrt{a_i})$, and the
rest of the arcs in a triangulation of $S \in R$ have length
comparable to $1$.
The period coordinates give the structure of a piecewise linear integral manifold on $\cQ\cM_g$. 
By work of Masur \cite{Masur1}, (up to normalization) the measure $\mu_{g}$ is the canonical measure defined by this piecewise linear 
integral structure (See \ch \ref{pc}).
Hence, by the definition of 
the measures $\mu_{g}$ and $\mm$ we get 
\begin{equation}
\label{eq:m:pi:R}
\mm(\pi(R)) \ge \mu(R) > c \prod_{i=1}^k a_i,
\end{equation}
where $c$ depends only on the genus. 
Note that by \cite{Rafi:TT} for any $Y \in \pi(R)$, the
only short simple closed curves on $Y$ (in the hyperbolic or extremal metric) 
are the core curves of the cylinders $C_i$, and the extremal
length of the core curve of $C_i$ is within a constant
multiple of $a_i$. Thus, there exists a constant $L'$ depending
only on the genus such that
\begin{displaymath}
\pi(R) \subset B_{\cT}(X,L').
\end{displaymath}
Note that the above equation takes place in $\cT_g/\Gamma_{g}$.
We may think of it as taking place in $\cT_g$ if we identify
$\pi(R)$ with a subset of the fundamental domain $\cF$. Then,
for any $\gg \in I_{X,L'}$, 
\begin{displaymath}
\gg \pi(R) \subset B_{\cT}(X,2L'). 
\end{displaymath}
Thus, in view of (\ref{eq:IXL}) and (\ref{eq:m:pi:R}), 
\begin{displaymath}
\mm(B_{\cT}(X,2L')) \ge |I_{X,L'}| \mm(\pi(R)) \ge c,
\end{displaymath}
where $c$ depends only on the genus. 
\qed\medskip
\subsection{Choosing a net in Teichm\"uller space} 
\label{netpoints}
 Let $c >2 L_{0}$ (the constant in Lemma ~\ref{lemma:uniform:vol:ball}),
$\widetilde{\cN}$ be a $(c,2c)-$net in $\cM_{g},$ and $\cN \subset \gp^{-1}(\widetilde{\cN})$ be a net in $\te_{g}$
such that $\gp(\cN)=\widetilde{\cN}.$
As before, $\gp: \te_{g} \rightarrow \cM_{g}$ is the natural projection to the moduli space. In other words, 
\begin{enumerate}
\item given $X \in \te_{g}$ there exists $Y \in \cN$
such that $d_{\cT}(X,Y) \leq 2 c$ , and 
\item for any $Y_{1}\not = Y_{2} \in \cN$, we have $d_{\cT}(Y_{1},Y_{2}) \geq c.$
\end{enumerate}
In view of Lemma~\ref{lemma:uniform:vol:ball}, for any $(c,2c)$
separated net $\cN$ in any Teichm\"uller space (including $\half^2$),
for $\tau \GG 1$ and any $X$, 
\begin{equation}
\label{eq:net:to:measure}
c_1 |B_{\cT}(X,\tau) \cap \cN| \le \mm(B_{\cT}(X,\tau)) \le c_2 |B_{\cT}(X,\tau) \cap
\cN|, 
\end{equation}
where $c_1, c_2$ depend only on $c$ and the genus. 
\begin{lemma}
Let $c> 2 L_{0}$. Then for any $(c,2c)$ net $\widetilde{\cN}$ in $\mathcal{M}_{g},$ there exists a constant $C_{2}>0$ such that for any $X
\in \te_{g}$ 
\begin{equation}\label{eq:netpoints}
 | \gp(B_{\cT}(X,\tau)) \cap \widetilde{\cN}| \leq C_{2} \tau^{3g-3}.
\end{equation}
\end{lemma}
\noindent
{\bf Proof.}
Fix $$\tau_{0}=\log\big(\frac{C_{g}^{2}}{\epsilon_{0}^{2}}\big)
$$ (see \ch \ref{boundedpants}). \\

\noindent
{\bf Step 1.} 
First we assume that $\tau\geq \tau_{0}.$
Let $\mathcal{A}_{X}=\{\alpha_{1},\alpha_{2},\ldots, \alpha_{k}\}$ be the set of simple closed curves of extremal length $\leq \epsilon_{0}^{2}\; e^{-2 \tau}$ on $X \in \cT_{g}.$ Let $\mathcal{P}_{1},\ldots,\mathcal{P}_{m}$ be the set of all combinatorially distinct pants decompositions of $S_{g}$ containing $\alpha_{1},\ldots, \alpha_{k}$. 

Given a pants decomposition $\mathcal{P}_{j}=\{\alpha_{1},\ldots, \alpha_{k},\alpha_{k+1}^{j},\ldots,\alpha_{3g-3}^{j} \} $ and ${\bf i}=(i_{1},\ldots, i_{3g-3})\in {\Bbb Z}^{3g-3}$ choose $Z_{j, \bf{i}} \in \cT_{g}$ such that: 
\begin{itemize}
\item 
for each $1 \leq l \leq k $
$$ e^{i_{l}-1} \leq \frac{\e_{\alpha_{l}}(Z_{j,\bf i})}{\e_{\alpha_{l}}(X)}\leq e^{i_{l}},$$
\item 
for $k+1 \leq l \leq 3g-3$ 
$$ C_{g} e^{i_{l}-1} \leq \e_{\alpha^{j}_{l}}(Z_{j,\bf i}) \leq C_{g} e^{i_{l}}.$$
\end{itemize}
Let
$$ \mathcal{I}= \{(i_{1},\ldots, i_{3g-3}) \in {\Bbb Z}^{3g-3}| \; \mbox{for}\;1 \leq  l \leq k, \; -2 \tau \leq i_{l}\leq 2\tau\; ,\mbox{and for}\; k+1\leq l \leq 3g-3\; , -5\tau+1\leq i_{l}\leq 0\ \}. $$
Define the set $\mathcal{Z}$ (with $|\mathcal{Z}|=O(\tau^{3g-3})$)  by  
$$\mathcal{Z}=\{Z_{j, (i_{1},i_{2},\ldots,i_{3g-3})} \;|\; {\bf i} \in \mathcal{I}\} \subset \cT_{g}.$$
\noindent
{\it Claim.} Given $Y \in B_{\cT}(X,\tau),$ there exists $Z_{j,{\bf i}}\in \mathcal{Z}$ such that 
$ \gp(Z_{j,{\bf i}}) \in \gp(B_{\cT}(Y, C))$, where $C$ is a constant which only depends on $g$.  This is a simple corollary of Theorem \ref{theorem:minsky:ext} and Theorem 
\ref{ker}.
Note that given 
$Y \in B_{\cT}(X,\tau)$ and $\alpha \in \mathcal{A}_{X}$, $\alpha$ is a short simple closed curve on $Y$ and  
 $$ e^{-2\tau} \leq \frac{\e_{\alpha}(Y)}{\e_{\alpha}(X)}\leq e^{2\tau}.$$
 Moreover, for any other simple close curve $\beta$ in a bounded pants decomposition of $Y$ since $d_{\cT}(X,Y) < \tau, $ if $\beta \not \in \mathcal{A}_{X}$ then 
 $$  \epsilon_{0}^{2} e^{-4 \tau} \leq \e_{\beta}(Y) \leq C_{g}. $$
 Now since any bounded pants decomposition $\mathcal{P}_{Y}$ on $Y$ has the combinatorial type of $\mathcal{P}_{j}$ for some $1\leq j \leq m,$ the claim follows easily from 
 Theorem \ref{theorem:minsky:ext}.\hfill$\Box$\\
 
As a result,  
\begin{equation}\label{O1}
 \gp(B_{\cT}(X,\tau)) \cap \widetilde{\cN} \subset \bigcup_{Z_{i} \in \mathcal{Z}} ( \gp(B_{\cT}(Z_{i},C)) \cap \widetilde{\cN}).
 \end{equation}
  \noindent
 {\bf Step 2.}
Since $\cN \subset \cT_{g}$ is a $(c,2c)$ net with $c > 2 L_{0}$, by Lemma $\ref{lemma:uniform:vol:ball}$ for any $Z \in \cT_{g}$
$$ | B_{\cT}(Z,C)) \cap \cN| \leq \frac{\vo(B_{\cT}(Z,C+c))}{\vo(B_{\cT}(Z,c))}= O(1),$$
and hence
  \begin{equation}\label{O2}
  |\gp(B_{\cT}(Z,C)) \cap \widetilde{\cN}|=O(1).
  \end{equation}
 Now since $|\mathcal{Z}|=O(\tau^{3g-3})$ the result follows from $(\ref{O1})$ and $(\ref{O2}).$

Note that $\tau_{0}$ only depends on $g$, and therefore if $\tau \leq \tau_{0}$ the result follows from step $2.$
\hfill $\Box$
 
\section{Geodesics in the thin part of moduli space}
\label{sec:thin:part}

In this section we prove Theorem~\ref{theorem:having:j:short:curves}.
The main idea, which is due to Margulis, is to prove a system of
inequalities, which shows that the flow (or more precisely an
associated random walk) is biased toward the thick part of
Teichm\"uller space. Variations on this theme have been used in
\cite{EMM}, \cite{EM} and \cite{Jayadev}.

\subsection{A system of inequalities}
\label{sec:ineq}

Suppose $0 < s < 1$ (in fact we will be using $s=1/2$ only).
Let $\tau \GG 1$ be a parameter to be chosen later. (In 
particular we will assume $e^{-(1-s)\tau} < 1/2$.)
Let $A_\tau$ be the
operator of averaging over a ball of radius $\tau$ in $\cT_{g}$ with respect to the Teichm\"uller metric. 
So if $f$ is a real-valued function on Teichm\"uller space,
then 
\begin{displaymath}
(A_\tau f)(X) = \frac{1}{\mm(B_{\cT}(X,\tau))} \int_{B_{\cT}(X,\tau)} f(Y) \,
d\mm(Y).
\end{displaymath}

\bold{Remark.} In \cite{EM}, \cite{Jayadev} and \cite{EMM} 
the average is over spheres. In this context, we use
balls, since Minsky's product region theorem gives us much more
precise information about balls than about spheres. 
\medskip

Let $m$ be the maximal number of disjoint simple closed curves on a closed surface
of genus $g$. Choose $\K > e^{2 m \tau}$, and pick constants
$\epsilon_1 < \epsilon_2 < \dots < \epsilon_m< 1/\K^3$ such that 
for all $1 \le i \le m-1$, 
\begin{equation}
\label{eq:choice:of:epsilons}
\epsilon_i < \frac{\epsilon_{i+1}}{\K^2 (2 m \K^2)^{2/s}}. 
\end{equation}
Note that $\K$ and $\epsilon_1, \dots, \epsilon_m$ are constants
which depend only on $\tau$ and the genus $g$. 

For $1 \le i \le m$ and $X \in \cT_g$ 
let $\ee_i(X)$ denote the extremal length of the $i$'th shortest simple closed 
curve on $X$. 
Let $f_0 = 1$ and for $1 \le j \le m$ let 
$$f_j(X) = \prod_{1 \le i \le j}
\left(\frac{\epsilon_i}{\ee_i(X)}\right)^s.$$  
Note that $f_j$ is invariant under the action of the mapping class
group, and thus descends to a function on $\cM_g$. 
Let
\begin{equation}\label{uu}
u(X) = \sum_{k=1}^m f_j(X). 
\end{equation}
Let $\epsilon_j' = \epsilon_j/(m \K^2)$. 
Let 
$$W_j = \{ X \in \cT_g \st \ee_{j+1}(X) > \epsilon_j' \}.$$

Note that $W_0$ is compact, and on $W_j$ there are at most $j$ short simple closed 
curves. If $X \not\in W_{j-1}$ then $X$ has at least $j$ short simple closed curves,
and thus if $X \in W_j \setminus W_{j-1}$ then $X$ has exactly $j$
short simple closed curves.

In this subsection, we prove the following:
\begin{proposition}
\label{prop:inequalities}
Set $s=1/2$. Then we may write 
\begin{equation}
\label{eq:main:inequality}
(A_\tau u)(X) \le c(X) u(X) + b(X), 
\end{equation}
where $b(X)$ is a bounded function which vanishes outside the compact
set $W_0$, and for all $j$ and for all $X \not\in W_{j-1}$, 
\begin{displaymath}
c(X)  \le C_j' \tau^j e^{-j \tau},
\end{displaymath}
where $C_j'$ depends only on the genus. 
\end{proposition}

We now begin the proof of Proposition~\ref{prop:inequalities}. If $\alpha_1, \dots, \alpha_j$ are disjoint simple closed curves 
let $P(\alpha_1,\dots, \alpha_j)$ denote the product region where for
all $1 \le i \le j$, the extremal length of $\alpha_i$ is at most $\K
\epsilon_i$.  

\begin{lemma}
\label{lemma:in:one:product:region}
Suppose $1 \le j \le m$, $1/2 \le s < 1$, and 
suppose $X$ and $\tau$ are such that $B_{\cT}(X,\tau)$ is completely
contained in $P(\alpha_1,\dots, \alpha_j)$. Also suppose
that for all $Y \in
B_{\cT}(X,\tau)$, the set $\{\alpha_1, \dots, \alpha_j\}$ coincides with the
set of the $j$ shortest curves on $y$. 
Then
\begin{displaymath}
(A_\tau f_j)(X) \le c_j f_j(X),
\end{displaymath}
where for $1/2 < s < 1$,
\begin{equation}
\label{eq:cj:bound}
c_j = C_j(s) e^{-2j(1-s) \tau}. 
\end{equation}
and if $s=1/2$, $c_j= C_j \tau^{j} e^{-j \tau}$. 
\end{lemma}

\bold{Sketch of proof.} 
%HERE AAAAaaaaaaaaaaaaaa
We will use Minsky's product region theorem as stated in Theorem $\ref{PRT}$. Let $\mm'$ denote the product measure
on $(\half^2)^j \cross \cT'$, and let $A'_\tau$
denote the averaging operator with respect to the product measure
$\mm'$, i.e. for a real-valued function $f$, 
\begin{displaymath}
(A'_\tau f)(X) = \frac{1}{\mm'(B_{\cT}(X,\tau))} \int_{B_{\cT}(X,\tau)} f(Y) \,
d\mm'(Y).
\end{displaymath}
We first establish the lemma with $A_\tau$ replaced by $A'_\tau$. 

We may write $X = (X_1,\dots,X_j,X')$ where $X_k$ is in the $k$'th
copy of the hyperbolic plane and $X' \in \cT'$. 
Because of the product region theorem 
$B_{\cT}(X,\tau)$ is essentially $B_{\cT}(X_1,\tau) \cross \dots \cross
B_{\cT}(X_j,\tau)  \cross B'(X',\tau)$ where for $1 \le k \le j$,
$B_{\cT}(X_k,\tau)$ is a ball of radius $\tau$ in the hyperbolic plane and
$B'(X,\tau)$ a ball of radius $\tau$ in $\cT'$. Also, 
by assumption, for any $Y \in B(X,\tau)$ the set of 
$j$ shortest simple closed curves on $Y$ is $\{\alpha_1,\dots,\alpha_j\}$. Thus,
for $Y \in B(X,\tau)$, with $Y = (Y_1,\dots, Y_j,Y')$, we have
\begin{displaymath}
f_j(Y) \approx (\epsilon_1 \dots \epsilon_j)^{s} 
\prod_{k=1}^j \ell_{min}(Y_k)^{-2s}
\end{displaymath}
where for $Y_{k} \in \half^2$, $\ell_{min}(Y_{k})$ is the flat length
of the shortest simple closed curve in the torus parametrized by $Y_{k}$.
We remark that the exponent is $2s$ instead of $s$ because on the torus
extremal length is the square of flat length. 
As a result, we have
\begin{equation}
\label{eq:A:prime:tau:fj}
(A'_\tau f_j)(X) \approx (\epsilon_1 \dots \epsilon_j)^s \prod_{k=1}^j
\left( \frac{1}{\Vol(B(X_k,\tau)}  \int_{B(X_k,\tau)} \ell_{min}(Y)^{-2s}
  \, d\Vol(Y) \right),
\end{equation}
where $\Vol$ is the standard volume form on $\half^2$.
 Now the integral in the parenthesis, 
i.e. an average of $\ell^{-2s}$ over a ball in
a hyperbolic plane is 
essentially done in \cite[Lemma 7.4]{EM} (except that there the
average is over spheres, but to get the average over balls one just
makes an extra integral over the radius). One gets for $1/2 < s < 1$, 
\begin{displaymath}
\frac{1}{\Vol(B(X,\tau)} 
\int_{B(X,\tau)} \ell_{min}(Y)^{-2s} \, d\Vol(Y) \le c(s) e^{-2(1-s)\tau}
\ell_{min}(X)
\end{displaymath}
and for $s=1/2$,
\begin{displaymath}
\frac{1}{\Vol(B(X,\tau)} 
\int_{B(X,\tau)} \ell_{min}(Y)^{-1} \, d\Vol(Y) \le c'\tau e^{\tau} \ell_{min}(X). 
\end{displaymath}
Substituting these expressions into (\ref{eq:A:prime:tau:fj})
completes the proof of Lemma~\ref{lemma:in:one:product:region}
with $A'_\tau$ instead of $A_\tau$.\\
In view
of the form of the function $f_j$, 
\begin{equation}
\label{eq:net:sum:to:integral}
c_1 \sum_{B(X,\tau) \cap \cN} f_j(Y) \le \int_{B(X,\tau)} f_j(Y) \,
d\mm(Y) \le c_2 \sum_{B(X,\tau) \cap \cN} f_j(Y).
\end{equation}
%HERE ADD STH.
Choose $L \GG \nb$ (where $\nb$ is as in
Theorem~\ref{PRT}), 
and choose a $(L,2L)$-separated net $\cN_k$
in each factor. Let $\cN$ be the product of the $\cN_k$.
Then $\cN$ is an $(L-\nb,2L+\nb)$-separated net in 
$P(\alpha_1,\dots,\alpha_j)$. Now in view of
(\ref{eq:net:to:measure}),  
\begin{displaymath}
\mm(B(X,\tau)) \approx |B(X,\tau) \cap \cN| \approx \prod_{k=1}^j 
|B(X_k,\tau) \cap \cN_k| \approx \prod_{k=1}^j \mm'(B(X_k,\tau))
\approx \mm'(B(X,\tau)),
\end{displaymath}
where, as before, $A\approx B$ means that the ratio $A/B$ is
bounded by two constants depending only on $\nb$, $L$ and $g$, and
thus ultimately only on $g$. Similarly, using
(\ref{eq:net:sum:to:integral}), we can show that
\begin{displaymath}
\int_{B(X,\tau)} f_j(Y) \, d\mm(Y) \approx \int_{B(X,\tau)} f_j(Y) \,
d\mm'(Y). 
\end{displaymath}
Thus $(A'_\tau f_j)(X) \approx (A_\tau f_j)(X)$.
\qed\medskip

\bold{Remark.} The proof works even if at some point $Y \in B(X,\tau)$
there are short simple closed curves other then $\{\alpha_1, \dots, \alpha_j\}$ (but
these other curves are longer then the maximum of the lengths of the
$\alpha_j$ at $Y$). This is used in the next lemma.  \medskip

\begin{lemma}
\label{lemma:inequalities}
For $1 \le j \le m$, let $u_j(X) = \sum_{k =j}^m f_j(X)$. 
Suppose $\ee_j(X) < \epsilon_j$. Then (assuming $\tau$ is large
enough), 
\begin{equation}
\label{eq:ineq:Atau:uj}
(A_\tau u_j)(X) \le \left(c_j+ \frac{1}{2\K}\right) u_j(X),
\end{equation}
where $c_j$ are as in (\ref{eq:cj:bound}). 
In particular, letting $j=1$, and noting that 
the set $\{ X \in \cM_g \st \ee_1(X) > \epsilon_1
  \}$ is compact, we have for all $X \in \cT_g$, 
\begin{displaymath}
(A_\tau u_1)(X) \le \left(c_1+ \frac{1}{2\K}\right)u_1(X) + b(\tau).
\end{displaymath}
\end{lemma}

\bold{Proof.} Note that for any $1 \le i \le m$, and any $X \in
\cT_g$, 
\begin{equation}
\label{eq:Atau:crude:bound}
\frac{1}{\K} f_i(X) \le (A_\tau f_i)(X) \le \K f_i(X)
\end{equation}
(this is because in $B(X,\tau)$ 
the extremal length of any simple closed curve cannot change by 
more than $e^{2 \tau}$). 

We divide the set $\{j, j+1, \dots m\}$ into two disjoint subsets: 
Let $I_1$ be the set of $k \in \{j, j+1, \dots m\}$ such that
\begin{equation}
\label{eq:insignificant:k}
f_k(X) \le \frac{u_j(X)}{2m\K^2},
\end{equation}
and let $I_2$ be the set of $k \in \{j, j+1, \dots m\}$ such that the
opposite inequality to (\ref{eq:insignificant:k}) holds. Suppose $k \in
I_1$. Then, by (\ref{eq:Atau:crude:bound}), $(A_\tau f_k)(X) \le
\frac{1}{2 m \K}f_k(X)$, and thus
\begin{equation}
\label{eq:contrib:insignificant}
\sum_{k \in I_1} (A_\tau f_k)(X) \le \frac{1}{2 m \K} u_j(X). 
\end{equation}

Now suppose $k \in I_2$. We claim that 
\begin{equation}
\label{eq:ellk:short}
\ee_k(X) < (2 m \K^2)^{1/s} \epsilon_k. 
\end{equation}
Indeed, if $k=j$ then (\ref{eq:ellk:short}) is true by assumption. 
If $k > j$ then 
\begin{displaymath}
u_j(X) \ge f_{k-1}(X) = f_k(X) \left(\frac{\ee_k(X)}{\epsilon_k}\right)^s \ge
\frac{u_j(X)}{2 m \K^2}  \left(\frac{\ee_k(X)}{\epsilon_k}\right)^s, 
\end{displaymath}
where we have used the inequality opposite to
(\ref{eq:insignificant:k}) in the last estimate. Thus
(\ref{eq:ellk:short}) follows. 

We now claim that under the assumption that $k \in I_2$ we have
\begin{equation}
\label{eq:ellk:plus:one:much:bigger}
\ee_{k+1}(X) \ge \K^2 \ee_k(X).
\end{equation}
If $k=m$ this is clear from (\ref{eq:ellk:short}) (since in the case
where $\ee_m(X)$ is small, there are no other short simple closed curves on $X$).
Now if $k < m$, then
\begin{displaymath}
u_j(X) \ge f_{k+1}(X) = f_k(X)
\left(\frac{\epsilon_{k+1}}{\ee_{k+1}(X)}\right)^s 
\ge  \frac{u_j(X)}{2 m \K^2} \left(\frac{\epsilon_{k+1}}{\ee_{k+1}(X)}\right)^s
\end{displaymath}
where again we used the inequality opposite to
(\ref{eq:insignificant:k}) in the last estimate. Thus,
\begin{equation}
\label{eq:ellk:plus:one:temp:long}
\ee_{k+1}(X) \ge \frac{\epsilon_{k+1}}{(2 m \K^2)^{1/s}}. 
\end{equation}
Now (\ref{eq:ellk:plus:one:much:bigger}) follows from
(\ref{eq:ellk:plus:one:temp:long}), (\ref{eq:ellk:short}) and 
(\ref{eq:choice:of:epsilons}). 

Now in view of (\ref{eq:ellk:short}) and
(\ref{eq:ellk:plus:one:much:bigger}), 
Lemma~\ref{lemma:in:one:product:region} can be applied to $f_k$. Thus,
for $k \in I_2$, 
\begin{equation}
\label{eq:contrib:significant:fk}
(A_\tau f_k)(X) \le c_k f_k(X) \le c_j f_k(X)
\end{equation}
where for the last inequality we assumed that $\tau$ was large
enough so that $c_k < c_j$ for $k > j$. 
Now (\ref{eq:ineq:Atau:uj}) follows from
(\ref{eq:contrib:insignificant}) and
(\ref{eq:contrib:significant:fk}). 
\qed\medskip

\bold{Proof of Proposition~\ref{prop:inequalities}.} 
Suppose $X \not\in W_{j-1}$. We may write
\begin{displaymath}
u(X) = u_j(X) + \sum_{k=1}^{j-1} f_k(X).
\end{displaymath}
Note that for $X \not\in W_{j-1}$ and $k < j$, 
$f_k(X) \le u_j(X)/(m \K^2)$. Hence, (\ref{eq:main:inequality}) follows
from (\ref{eq:ineq:Atau:uj}) and (\ref{eq:Atau:crude:bound}). 
\qed\medskip

\subsection{A uniform estimate for the measure of a ball}
\label{sec:measures:balls}
Given $X$ and $Y$ in $\cT_g$ and $\tau > 0$, let 
$$F_\tau(X,Y)= \{ Z\; |\;  Z \in \Gamma_{g} \cdot Y, \; d_{\cT}(X,Z) \leq \tau \} \subset \cT_{g}.$$ In the Appendix we show:
\begin{prop}
%\label{prop:toprove}
Given $\epsilon>0,$ there exists $\tau_{0}>0$ such that for any $\tau>\tau_{0}$ and $X,Y \in \te_{g}$
we have 
$$ | F_{\tau}(X,Y)| \leq e^{(h+\epsilon) \tau} G(Y)^{2}.$$
\end{prop}
In order to prove this statement we show that in general
\begin{equation}\label{sv}
 |F_{\tau}(X,Y)|=O( \tau^{3g-3} e^{h \tau} G(X) G(Y)).
 \end{equation}
Using these estimates we show that:
\begin{proposition}
\label{prop:measures:of:balls}
There exists a constant $C_2$ such
that for any $X$, any $\delta'' > 0$ 
and any sufficiently large $\tau$, the volume of any
$B(X,\tau)$ is bounded by $C_2 e^{(h+\delta'')\tau}$.  
\end{proposition}
\noindent
\bold{Proof.} See Appendix.
\qed\medskip\\
\noindent
{\bf Remark.} Lemma $\ref{theorem:uniform:upper:bound}$ in the appendix is a stronger version of $(\ref{sv})$. The proof of this lemma uses Theorem $\ref{theorem:minsky:ext}$. This statement can be generalized for other strata of moduli space of quadratic differentials. 

\subsection{Proof of Theorem~\ref{theorem:having:j:short:curves}} 
\label{sec:typeII}
As in \ref{netpoints}, we discretize Teichm\"uller space by fixing a $ (c,2c)$
separated net $\cN \subset \cT_g$.
We note that there exists a constant $\kappa_0$ such that
for all $X \in \cT_g$, 
$$
\frac{1}{\kappa_0}\; u(X) \le G(X) \le \kappa_0\; u(X),
$$
where $\kappa_0$ depends on $\tau$ and the $\epsilon_i$ (and thus 
ultimately only on $\tau$ and the genus). Also, as in (\ref{GG}),
\begin{displaymath}
G(X) = \prod_{1 \le i \le m} \ee_i(X)^{-1/2}. 
\end{displaymath}

\bold{Trajectories of the random walk.}
Suppose $R \GG \tau$ and let $n$ be the integer part of $R/\tau$. 
By a {\it trajectory of the random walk} we mean a map $\lambda: \{ 0, n-1\}
\to \cT_g$ such that :
\begin{itemize}
\item for all $0 < k \le n-1$ we have $d(\lambda_k,
\lambda_{k-1}) \le \tau$, and 
\item $\lambda_k$ belongs to the net $\cN$.
\end{itemize}
Let $\cP_{\tau}(X,R)$ denote the set of all trajectories for which
$d(\lambda_0,X) = O(1)$. It is a corollary of
Proposition~\ref{prop:measures:of:balls} that
\begin{equation}
\label{eq:number:all:trajectories}
|\cP_{\tau}(X,R)| \le C_2 e^{(h+\delta'')R},
\end{equation}
where $| \cdot|$ denotes the cardinality of a set. 

We say that a trajectory is {\em almost closed in the quotient} if the
distance in $\cM_g$ between the projection to $\cM_g$ of $\lambda(0)$
and the projection to $\cM_g$ of $\lambda_{n-1}$ is $O(1)$. 

Let $\delta > 0$ be a constant to be chosen later. (We will have
$\delta < \epsilon_j'$ for $1 \le j \le m$, where the $\epsilon_j'$
are as in \S\ref{sec:ineq}). For $j \in 
\natls$, let
$\cP_{j,\tau}(X,\delta,R)$ 
denote the set of all trajectories starting within $O(1)$
of $X$ for which at any point, there are at least $j$ simple closed curves of length
at most $\delta$. Let $\widetilde{\cP}_j(X,\delta,R)$ denote
the subset of these trajectories which are almost closed in the
quotient. 
 We use the systems of inequalities discussed in $\ch \ref{sec:ineq}$ to show:

\begin{lemma}
\label{lemma:est:tilde:P:j} 
Let $j \in \natls$. For any $\epsilon' >
0$, there exists $\tau_{0}>0$ such that for $\tau> \tau_{0}$, $\delta
> 0$ small enough (depending on $\tau$, $\epsilon'$ and $g$),  and any sufficiently large $R$ (depending only on $\epsilon'$ and $\tau$)
we have:
\begin{equation}
\label{eq:est:tilde:P:j}
|\widetilde{\cP}_{j,\tau}(X,\delta,R)| \le e^{(h-j+\epsilon') R}.
\end{equation}
% WHAT C EXISTS IN THE BEGINNING//\\Here $C$ depends on $\cN$, $\delta$, $\epsilon'$ and $g$.
\end{lemma}
%C is indep. of \tau emphesis..............................................
\bold{Proof.} 
By the definition of $u$ for $s=1/2$ (see \ch \ref{sec:ineq}), it is easy to check that 
for any $Y \in \cT_{g}$
\begin{equation}
\label{eq:G:same:as:u}
\frac{1}{\kappa_0} \; u(Y) \leq G(Y) \leq \kappa_0 \; u(Y),
\end{equation}
where $$\kappa_0 \leq e^{m_{0} \tau}.$$ Here $m_{0}$ is a universal constant independent of $\tau$ and $Y$. 
Also, by the definition of the function $u$ if $d_{\cT}(X_{1},X_{2})=O(1),$ then 
\begin{equation}\label{cl}
\frac{u(X_{1})}{u(X_{2})}=O(1).
\end{equation}
Note that if $\lambda \in \widetilde{\cP}_j(X,\delta,R)$ then
$G(X) \approx G(\lambda_{n-1}) \le \kappa_0 u(\lambda_{n-1})$. 
Let $R = n \tau$, and let
\begin{displaymath}
q_j(X,R) = \sum_{\lambda \in \cP_j(X,\delta,R)} u(\lambda_{n-1})
\end{displaymath}
Therefore,
\begin{equation}
\label{eq:frac:estimate}
|\widetilde{\cP}_j(X,\delta,R)| \le 
\frac{C_1}{G(X)} 
\sum_{\lambda \in \widetilde{\cP}_j(X,\delta,R)} \kappa_0 u(\lambda_{n-1})
= \frac{C_1 \kappa_0}{G(X)} q_j(X,R),
\end{equation}
where $C_1 = O(1)$. For $0 < r=k\tau < R$ let 
\begin{displaymath}
q_j(X,R,r) = \sum_{\lambda \in
  \cP_j(X,\delta,R,r)} u(\lambda_{k-1}), 
\end{displaymath}
where the elements of $\cP_j(X,\delta,R,r)$ 
are the trajectories $\lambda$
belonging to $\cP_j(X,\delta,R)$ but truncated after 
$k = r/\tau$ steps.
Then
\begin{equation}\label{b}
q_j(X,R,r+\tau)  =  \sum_{\lambda \in
  \cP_j(X,\delta,R,r+\tau)} u(\lambda_k) 
  \end{equation}
\begin{align*}
= \sum_{\lambda \in \cP_j(X,\delta,R,r)}
 \qquad \sum_{\lambda_k \in \cN \cap
  B(\lambda_{k-1},\tau)} u(\lambda_k) \\
& \le C \sum_{\lambda \in \cP_j(X,\delta,R,r)}
\int_{B(\lambda_{k-1},\tau)} u(Y) \, d\mm(y) \\
& = C \sum_{\lambda \in \cP_j(X,\delta,R,r)} \mm(B(\lambda_k,\tau))
\ (A_\tau u)(\lambda_{k-1})  
\end{align*}
where in the next to last line we use (\ref{cl}) and Lemma \ref{lemma:uniform:vol:ball} to estimated a sum over $\cN \cap
B(\lambda_{k-1},\tau)$ by a constant $C$
times an integral over $B(\lambda_{k-1},\tau)$.

Note that for $\lambda \in \cP_j(X,\delta,R)$, 
the number of simple closed curves shorter then $\delta$ 
on $\lambda_{k-1}$ is at least $j$. Thus, if $\delta$ is small enough,
(depending on the the $\epsilon_j'$ and thus ultimately only on
$\tau$ and the genus), $\lambda_{k-1} \not\in W_{j-1}$. 
Then, from Proposition~\ref{prop:inequalities}, and assuming $\tau$ is
large enough so that Proposition~\ref{prop:measures:of:balls} holds 
with $\delta'' < \epsilon'/3$, we get 
\begin{equation}
\label{eq:q:Delta:iter}
q_j(X,R,r+\tau) \le C C_{2} \sum_{\lambda \in
  \cP_j(X,\delta,R,r)} e^{(h+\epsilon'/3)\tau} C_j'
\tau e^{-j\tau} u(\lambda_{k-1}) = C C_j' C_{2}
\tau e^{(h+(\epsilon'/3)-j)\tau} q_j(X,R,r).
\end{equation}
Now iterating (\ref{eq:q:Delta:iter}) $n= R/\tau$ times we get
\begin{equation}
\label{eq:qj:estimate}
q_j(X,R) \le u(X) (C C_{2}\;C_j')^n \tau^{n} e^{(h-j+\epsilon'/3) n\tau} 
=  u(X) e^{(h -j + \epsilon'/3 + \frac{\log
  (C C_{2} C_j'\tau)}{\tau}) R}. 
\end{equation}
Hence, in view of (\ref{eq:frac:estimate}) and (\ref{eq:G:same:as:u}), we get 
\begin{equation}\label{final}
|\widetilde{\cP}_j(X,\delta,R)| \le   C_{1} e^{2 m_{0} \tau} \times e^{(h -j + \epsilon'/3+ \frac{\log
  (C C_{2} C_j'\tau)}{\tau}) R} \leq  C_{1} e^{(h -j + \epsilon'/3 + \frac{\log (C C_{2} C_j'\tau)}{\tau}+ 2m_{0}\; \frac{\tau}{R} ) R}.
  \end{equation}
Now it is enough choose $\tau_{0}$ so that $\frac{\log (C\; C_{2}\; C_j'\tau_{0})}{\tau_{0}} < \epsilon'/3$. Then for any $\tau \geq \tau_{0}$ and $R \geq \frac{6m_{0} \tau+ 3 \log(C_{1})}{\epsilon'}  $
(\ref{final}) implies $(\ref{eq:est:tilde:P:j}).$ 

\qed\medskip\\

\noindent
{\bf Remark.}
One can use Lemma $\ref{theorem:uniform:upper:bound}$ to prove Lemma $\ref{lemma:est:tilde:P:j}$.

Let $N_j(X,\delta,R)$ be the number of conjugacy classes of closed 
geodesics of length at most $R$ which pass within $O(1)$ of the point
$X$ and always have at least $j$ simple closed curves of length at most $\delta$.
%CHHHHHHHHHHHHHHHHHHHHHHHHHHHHHHHHHHHHHHHHHHHHOW ARE THESE DEPEND ON NET???????????
\begin{lemma}
\label{lemma:geodeiscs:to:trajectories}
For any $\epsilon' > 0$ we may choose $\tau$ large enough (depending
only on $\epsilon'$) so that for any $X \in \cT_g$, any $\delta < 1/c_{0}$
and any sufficiently large $R$ (depending only on $\epsilon'$ and $\tau$)
we have
\begin{equation}
\label{eq:Nx:tilde:S}
N_j(X,\delta,(1-\epsilon')R) < C |\widetilde{\cP}_j(X, c_{0} \delta,R)|,
\end{equation}
where $C$ and $c_{0}$ are constants which only depend on $\mathcal{N}$ and $g$. 
%CCCCCCCHHHHHHHEEEEEECCCCCCKKKKKKK\\
\end{lemma}

\bold{Proof.} Let $$I_X=\{ \gg \in \Gamma_{g}\; | B_{\cT}(X, \gg \cdot X) \leq 1\} .$$ In other words, $I_{X}$ is the 
subset of the mapping class group which moves $X$ by at most
$O(1)$. In fact, $| I_X| \approx  G(X)^2$ (See also 
Theorem $\ref{theorem:uniform:upper:bound}$).

% The number of Dehn twists around
% a curve $\beta \in S(x)$ which move a point $O(1)$ is
% $1/\e_{\beta}(x)$. Therefore, $|I_x| \approx G(x)$. 
Now consider a closed geodesic $\gamma$ in $\cM_g$ which passes
within $O(1)$ of $\gp(X)$ (recall that $\gp$ denotes the projection map 
$\gp: \cT_g \to \cM_g$). Let $[\gamma]$ denote the corresponding
conjugacy class in $\Gamma_{g}$. Then there are approximately $|I_X|$
lifts of $[\gamma]$ to $\cT_g$ which start within $O(1)$ of $X$. Each
lift $\tilde{\gamma}$ is a geodesic segment of length equal to the length of
$\gamma$.

We can mark points distance $\tau$
apart on $\tilde{\gamma}$, and replace these points
by the nearest net points. (This replacement is the cause of the
$\epsilon'$). This gives a map $\Psi$ 
from lifts of geodesics to trajectories.
If the original geodesic $\gamma$ has length at most $(1-\epsilon')R$
and always has $j$ simple closed curves shorter then $\delta$, 
then by Theorem \ref{ker} the resulting trajectory belongs
to $\widetilde{\cP}_j(X, c_{0} \delta,R)$, where $c_{0}$ only depends on $\mathcal{N}.$ 
%(CCCCCCCCCCHHHHHHHHEEEEEECCCCCCCKKKKKKKKKKKK)

If two geodesic segments map to the same trajectory, then the segments
fellow travel within $O(1)$ of each other. In particular if $\gg_1$ and
$\gg_2$ are the pseudo-Anosov elements corresponding to the two
geodesics, then $d_{\cT}(\gg_2^{-1} \gg_1 X, X) = O(1)$, thus $\gg_2^{-1} \gg_1 \in
I_X$.

We now consider all possible geodesics contributing to
$N_j(X,\delta,(1-\epsilon')R)$; for each of these we consider all the
possible lifts which pass near $X$, and then for each lift consider
the associated random walk trajectory. We get:
\begin{displaymath}
N_j(X,\delta,(1-\epsilon')R) |I_X| \le C |I_X|
|\widetilde{\cP}_j(X, c_{0}\delta,R)|
\end{displaymath}
(the factor of $|I_X|$ on the left hand side is due to the fact that
we are considering all possible lifts which pass near $X$, and the
factor of $|I_X|$ on the right is the maximum possible number of times
a given random walk trajectory can occur as a result of this process). Thus,
the factor of $|I_X|$ cancels, and the lemma follows.
\qed\medskip

The following lemma is due to Veech \cite{Veech}.
\begin{lemma}
\label{lemma:veech}
Suppose $\gamma$ is a closed geodesic of length at most
$R$ on $\cM_g$. Then for any $X \in \gamma$ the extremal length of the shortest simple closed curve on $X$ is at least $\epsilon_0' e^{-(6g-4)R}$, where $\epsilon_0'$ depends only on $g$. 
\end{lemma}

\bold{Proof.} We reproduce the proof for completeness. Let $\tilde{X}$
be some point in $\cT_g$ with $\gp(\tilde{X}) = X$. 
Suppose the estimate is false, and let $\alpha$ be a simple closed curve on
$\tilde{X}$ with extremal length less than $\epsilon_0' e^{-(6g-4)R}$. Let
$\gg$ be the element of the mapping class group associated to
the lift of $\gamma$ passing through $\tilde{X}$. 
$\e_\alpha(\tilde{X}) \le
\epsilon_0' e^{-(6g-4)R}$. Therefore, by Theorem $\ref{ker}$, for 
$j \in \natls$, 
\begin{displaymath}
\e_{\gg^j \alpha}(\tilde{X}) = \e_\alpha(\gg^{-j} \tilde{X}) \le
\epsilon_0' e^{-(6g-4+2j)R}
\end{displaymath}
In particular, for $1 \le j \le 3g -2$, $\e_{\gg^j
\alpha}(\tilde{X}) \le \epsilon_0'$. Therefore (assuming that
$\epsilon_0'$ is small enough), for $1 \le j\le 3g-2$,the simple closed curves 
$\gg^j \alpha$ are disjoint. This is a contradiction.
\qed\medskip\\
\bold{Remark.} 
Recall that by $(\ref{Maskit})$ for very short curves the hyperbolic and the extremal
lengths are comparable. Therefore hyperbolic length of the shortest closed geodesic on $X \in \gamma$
is at least $\epsilon_0'' e^{-(6g-6)R},$ where $\epsilon_0''$ depends only on $g$.
 
\bold{Proof of Theorem~\ref{theorem:having:j:short:curves}.}
Let $\epsilon' = \epsilon/8c_{0}$, where $c_{0}$ is the constant in Theorem $\ref{lemma:geodeiscs:to:trajectories}.$ 
By Lemma~\ref{lemma:geodeiscs:to:trajectories} and
Lemma~\ref{lemma:est:tilde:P:j} we can choose $\tau$ and $\delta$
so that (\ref{eq:Nx:tilde:S}) holds and also (\ref{eq:est:tilde:P:j})
holds with $\delta$ replaced by $c_{0} \delta$. We get, 
for sufficiently large $R$, 
\begin{equation}
\label{eq:Nx:Q}
N_j(X,\delta, R) < C' e^{(h-j+\epsilon' c_{0}/4)R}.
\end{equation}
%Is this OK here? because of the change from 2\delta to c_0 \delta-------??????????????????
Finally $N_j(\delta,R)$ is at most 
$\sum_X
N_j(X,\delta,R)$, where we have to let $X$ vary over the net points within distance $1$ of a
fundamental domain for the 
action of the mapping class group. In view of Lemma~\ref{lemma:veech},
the number of relevant points in the net is
at most polynomial in $R$. Thus
Theorem~\ref{theorem:having:j:short:curves}
follows.
\qed\medskip

\section{Recurrence of Geodesics}\label{rec}

In this section, we discuss basic recurrence properties of closed geodesics in 
$\cQ\cM_g.$
Note that a Teichm\"uller geodesic $\gamma$ is in fact a path in the the space of unit area holomophic quadratic differentials on surfaces of
genus $g$. As in \S \ref{qd}, let $\cQ\cM_{g}(1,\dots,1) \subset \cQ\cM_{g}$ denote the {\em principal stratum}, i.e. the set of pairs $(X,q)$ where $q$ is a holomoprhic
quadratic differential on $X$ with simple zeroes.  

\noindent
{\bf Notation.}
For a compact subset $\cKM$ of $\cM_g$ and a number $\theta > 0$ let
$N^{\cKM}(\theta,R)$ denote the number of closed geodesics $\gamma$ of
length at most $R$ such that $\gamma$ spends at least
$\theta$-fraction of the time outside $\cKM$. 

Similarly, for a compact subset $\cKQ$ of $\cQ\cM_{g}(1,\dots,1)$ and $\theta > 0$ we denote by 
$N^{\cKQ}(\theta,R)$ the number of closed geodesics $\gamma$ of
length at most $R$ such that $\gamma$ spends at least
$\theta$-fraction of the time outside $\cK$. 

We prove the following theorems:

\begin{theorem}
\label{theorem:staying:in:think:part}
Suppose $\theta > 0$. Then there exists a compact subset $\cKM$ of
$\cM_g$ and $\delta > 0$ such that for sufficiently large $R$, 
\begin{displaymath}
N^{\cKM}(\theta,R) \le e^{(h-\delta)R}.
\end{displaymath}
\end{theorem}

\begin{theorem}
\label{theorem:staying:in:principal:stratum}
Suppose $\theta > 0$. Then there exists a compact subset $\cKQ \subset
 \cQ\cM_{g}(1,\dots,1)$ and $\delta > 0$ such that for sufficiently large $R$, 
\begin{displaymath}
N^{\cKQ}(\theta,R) \le e^{(h-\delta)R}.
\end{displaymath}
\end{theorem}

\bold{Proof of \ref{theorem:staying:in:think:part}.} In view of
Corollary~\ref{cor:number:staying:outside:compact}, there exists
a compact set $\cKM_1$, such that the number of geodesics which
do not return to $\cKM_1$ is $O(e^{(h-0.99)R})$. But then using
Proposition~\ref{prop:inequalities} and \cite[Corollary 2.7]{Jayadev}
(cf. \cite[Theorem 2.3]{Jayadev})
\mc{give exact reference}, there exists a compact
set $\cKM$ depending on $\cKM_1$, $\theta$ and there exists $\delta' >
0$, such that the number of random walk trajectories which start in
$\cKM_1$ and spend at least
$\theta$-fraction of the time outside of $\cK$ is at most
$e^{(h-\delta')R}$. (Note that even though 
\cite[Corollary 2.7]{Jayadev} is not stated with an exponential bound,
the proof does in fact imply this, as is done in the statement and
proof of \cite[Theorem 2.3]{Jayadev}). 
It follows that the same kind of
estimate is true for the number of geodesics (see the proof
of Lemma~\ref{lemma:geodeiscs:to:trajectories}). 
\qed\medskip 

The rest of this subsection will consist of the proof of
Theorem~\ref{theorem:staying:in:principal:stratum}. 
\subsection {Choosing hyperbolic neighborhoods of points}\label{CNP}
\noindent
{\bf Notation.} Given $q\in \cQ\cM_{g}$ we let $\cl(q)$ denote the length of the shortest saddle connection on $q$
in the flat metric defined by $q$. Suppose $\cKM \subset \cM_g$ is a compact set.
 For simplicity, we denote the preimage of $\cKM$ in $\te_{g}$ by $\cKT=\cP^{-1}(\cKM)$.
Also, to simplify the notation, let $q_{X,Y} \in \mathcal{Q}^{1}(X)$ be the quadratic differential such that $g_{d_{\cT}(X,Y)}(q) \in \mathcal{Q}^{1}(Y).$

As in \cite[\S{2}]{ABEM}, we denote the strong unstable, unstable,
stable and strong stable foliations of the geodesic flow 
by $\cF^{uu}$, $\cF^{u}$,
$\cF^{s}$ and $\cF^{ss}$ respectively; for a 
given quadratic differential $q$, 
$$\cF^{ss}(q)=\{ q_{1} \in \mathcal{Q}^{1}\te_{g}\;|\; \RR(q_{1})=\RR(q)\},  \cF^{s}(q)=\{ q_{1} \in \mathcal{Q}^{1}\te_{g}\;|\; [\RR(q_{1})]=[\RR(q)]\},$$
and 
$$\cF^{uu}(q)=\{ q_{1} \in \mathcal{Q}^{1}\te_{g}\;|\; \II(q_{1})=\II(q)\}, \cF^{u}(q)=\{ q_{1} \in \mathcal{Q}^{1}\te_{g}\;|\; [\II(q_{1})]=[\II(q)]\}.$$
%PLEASE CHECK THIS SENTENCE??????????????????????????????????????????????????????????????????????????????????????????????????
We consider the distance function defined by the modified Hodge norm $d_{H}$ on each horosphere $\cF^{ss}$. This is closely related to the Hodge norm studied by 
Forni \cite{forni}.
We also consider $d_E( \cdot, \cdot)$ the Euclidean metric on each horosphere as defined in \cite[\S{3.5}]{ABEM}. In fact, the Euclidean norm on the tangent space of $\mathcal{Q}^1\mathcal{M}_{g}$ is defined using period coordinates (see $\ch \ref{pc}$).
We remark that this norm depends on the choice of a triangulation on the surface. However,  in a given compact subset of $\mathcal{Q}^1\mathcal{M}_{g}$
 it is well-defined up to a multiplicative constant. So we can use it to measure the "distance" between two quadratic differentials in a compact subset of 
 $\mathcal{Q}^1\mathcal{M}_{g}.$
 We will show:
\begin{lemma}
\label{lemma:compactness}
Suppose $\cKM \subset \cM_g$ is compact. Given $\rho,\epsilon >0,$ and $1> \nb>0$, there
exists $\rho_{0}>0$ (depending only on $\cKM$, and $\nb$), and
$L_{0}=L_{0}(\cKM,\epsilon,\rho,\nb)$ such that if $X, p_{0} \in \cKT$, $d_{\cT} (p_{0}, p_{1}) < \rho_{0},$
$d_{\cT}(X,p_{1})=L > L_{0}$, and
$$\{ s\; | s \in [0, L] \; ,\; \cl(g_{s}(q_{X,p_{0}})) \geq \epsilon \} > \nb L,$$
then $$d_{E}(q,q_{X,p_{0}}) < \rho,$$
where $q$ is the unique quadratic differential in $ \cF^{uu}(q_{X,p_{0}}) \cap \cF^{s}(q_{X,p{1}})$.
\end{lemma} 
\noindent
{\bf Remark}. Note that given $q_{1},q_{2} \in \mathcal{Q}^{1}\te_{g}$, $|\cF^{uu}(q_{1}) \cap \cF^{s}(q_{2})| \leq 1$. In general, this set can be empty, but if $q_{2}$ is close enough to $q_{1}$ then $\cF^{uu}(q_{1}) \cap \cF^{s}(q_{2})\not=\emptyset$. The proof of Lemma $\ref{lemma:compactness}$ implies that when $\rho_{0}$ is small and $L_{0}$ is large 
 $ \cF^{uu}(q_{X,p_{0}}) \cap \cF^{s}(q_{X,p_{1}}) \not =\emptyset.$\\
We will show:
\begin{lemma}
\label{lemma:hyperbolic:nbhd}
Suppose $\cKM \subset \cM_g$ is compact. Given $\epsilon > 0$, there exists a constants $L_0$, depending only on $\epsilon$ and
$\cKM$, and $c_{0}$ depending only on $\cKM$ with the following property.

If  \begin{itemize}
\item 
$\gamma: [0,L] \to \cQ\cT_g$ is a geodesic segment (parametrized by
arclength) such that $(\gp \circ \pi)(\gamma(0)) \in \cKM, (\gp \circ \pi)(\gamma(L)) \in \cKM $,
and $L> L_{0},$
% WHy was this K_1 :::::::::::::::::::??????????????? such that $(\gp \circ \pi)(\gamma(0)) \in \cKM_1, (\gp \circ \pi)(\gamma(L)) \in \cKM_1$,
\item  $\hat{\gamma}: [0,L'] \rightarrow \cQ\cT_{g}$ is the geodesic segment connecting $p_1, p_{2} \in \cT_{g}$ 
such that 
$$d_{\cT}(p_1,\pi(\gamma(0)) < c_0\;\; , \;\; d_{\cT}(p_2,\pi(\gamma(L))) < c_0,$$
and 
\item 
\begin{displaymath}
|\{ t \in [0,L] \st \cl(\gamma(t)) > \epsilon \}| > (1/2) L,
\end{displaymath}
\end{itemize}
then 
\begin{displaymath}
|\{ t \in [0,L'] \st \cl(\hat{\gamma}(t)) > \epsilon/4 \}| > (1/3) L'.
\end{displaymath}
\end{lemma}

\subsection{Hodge and Euclidean distance functions on the stable and
  unstable foliations}  
First, we briefly recall some useful decay properties of the Hodge and
Euclidean distances proved in \cite[\S{3.4}, \S{3.5}]{ABEM}:
\begin{itemize}
\item {\bf P1}: Given $\epsilon>0$, there exists $\epsilon'>0$ such that if $\cl(q_{1}) \geq \epsilon$ and $d_{E}(q_{1},q_{2}) \leq \epsilon'$ then 
$\cl(q_{2}) \geq \epsilon/2.$ 
\item {\bf P2:} Assume that $\cKT_{1} \subset \te_{g}$ is the preimage 
of a compact subset $\cKM_{1} \subset \cM_{g}.$ Given $\rho_{1}>0,$ there exists $\rho_{0}$ such that if $q_{1} \in \cF^{ss}(q_{2})$, $\pi(q_{1}) \in \cKT_{1},$ and 
$d_{E}(q_{1},q_{2})> \rho_{1},$ then $d_{\cT}(\pi(q_1),\pi(q_2))> \rho_{0}.$

\item {\bf P3:}
 There exists $C_{1}>0$ such that if $d_{H}(q_{1},q_{2})<1$,  $q_{1} \in \cF^{ss}(q_{2})$ and $s\geq 0$
\begin{equation}
\label{one} 
d_{E}( g_{s}q_{1}, g_{s}q_{2}) \leq d_{H}( g_{s}q_{1}, g_{s}q_{2}) < C_{1} d_{H} (q_{1},q_{2});
\end{equation}
\item {\bf P4:} Moreover, 
given $\epsilon,\nb>0$ , there exist $C_{0},C_{0}'$ and $\na>0$ such that for any $q_{1} \in \cF^{ss}(q_{2})$ with $d_{H}(q_{1},q_{2})<1$, and $s\geq 0$ 
if 
\begin{equation}
\label{con}
|\{ t | \; t \in [0,s] , \cl(g_{t}q_{1}) > \epsilon \}| > \nb\; s,
\end{equation}
 and $\pi(g_{s}q_{1}) \in \cKT_{1}$, then
\begin{equation}
\label{two}
d_{E}( g_{s}q_{1}, g_{s}q_{2})\leq C_{0}' d_{H}(g_s q_{1},g_s q_{2}) < C_{0} e^{-\na s} d_{H} (q_{1},q_{2}).
\end{equation}
Note that in this case, by $(\ref{two})$, there exists $L_{0}$ (depending only on $\cKT_{1},$
$\nb$ and $\epsilon$) such that for $s> L_{0}$, $(\ref{con})$ implies that :
\begin{equation}
|\{ t | \; t \in [0,s] , \cl(g_{t}q_{2}) > \epsilon/2 \}| > \frac{\nb}{2} s.
\end{equation}
\end{itemize}
\noindent
% Note that the flat metric associated to 
% any $q \in \cF^{u}_L$ near \mc{how near} $\gamma(L)$ does not contain
% any horizontal saddle connections of bounded length (or else $g_t^{-1} q$
% will have a small saddle connection for $t \GG 0$, contradicting the
% assumption ????. 
\noindent

\makefig{Proof of Lemma~\ref{lemma:compactness}.}{fig:projection}{\begin{picture}(0,0)%
\includegraphics{projection.pstex}%
\end{picture}%
\setlength{\unitlength}{3947sp}%
\begingroup\makeatletter\ifx\SetFigFont\undefined%
\gdef\SetFigFont#1#2#3#4#5{%
  \reset@font\fontsize{#1}{#2pt}%
  \fontfamily{#3}\fontseries{#4}\fontshape{#5}%
  \selectfont}%
\fi\endgroup%
\begin{picture}(7824,4154)(2089,-6444)
\put(4576,-3436){\makebox(0,0)[lb]{\smash{{\SetFigFont{12}{14.4}{\rmdefault}{\mddefault}{\updefault}{\color[rgb]{0,0,0}$U$}%
}}}}
\put(6001,-2461){\makebox(0,0)[lb]{\smash{{\SetFigFont{12}{14.4}{\rmdefault}{\mddefault}{\updefault}{\color[rgb]{0,0,0}$\gamma(0)$}%
}}}}
\put(3601,-4111){\makebox(0,0)[lb]{\smash{{\SetFigFont{12}{14.4}{\rmdefault}{\mddefault}{\updefault}{\color[rgb]{0,0,0}$q$}%
}}}}
\put(3226,-2536){\makebox(0,0)[lb]{\smash{{\SetFigFont{12}{14.4}{\rmdefault}{\mddefault}{\updefault}{\color[rgb]{0,0,0}$f(q)$}%
}}}}
\put(2701,-5761){\makebox(0,0)[lb]{\smash{{\SetFigFont{12}{14.4}{\rmdefault}{\mddefault}{\updefault}{\color[rgb]{0,0,0}$S_0$}%
}}}}
\put(9151,-2611){\makebox(0,0)[lb]{\smash{{\SetFigFont{12}{14.4}{\rmdefault}{\mddefault}{\updefault}{\color[rgb]{0,0,0}$\cF_0$}%
}}}}
\end{picture}%
}

\bold{Proof of Lemma $\ref{lemma:compactness}$.} 
Let $B_E(q,r)$ denote the ball of radius $r$ with center $q$ in the
Euclidean metric. 
%Since $\cKM_1$ is compact, there is a number
%$\rho_1 > 0$ depending only on $\cKM_1$ and $\epsilon$ such that 
%$B_E(q,\rho_1)$ is contained in one
%fundamental domain for the action of the mapping class group
%$\Gamma$.
Let $S_0 = S(X)$ denote 
the sphere at $X$, i.e. the set of unit area
holomorhic quadratic differentials on the surface $X$, and $\gamma(t)=g_{t}(q_{X,p_{0}}).$
For $q \in S_0$ near $\gamma(0)=q_{X,p_{0}}$, 
let 
$$f(q)=  \cF^{uu}(q_{X,p_{0}}) \cap \cF^{s}(q);$$ 
in other words, we can choose $t(q) \in \reals$ be such that
$f(q) \in \cF^{uu}(q_{X,p_{0}})$ 
and $g_{t(q)} q$ and $f(q)$ are on the same leaf of $\cF^{ss}$.
Then clearly $f(q_{X,p_0}) = q_{X,p_0}$, and there exists a number $\rho_3 > 0$ depending only on $\cKM$ \mc{check
  this again. worry about branches} such that
the restriction of $f$ to $B_E(q_{X,p_0},\rho_3) \cap S_0$ is a
homeomorphism onto a neighborhood of $q_{X,p_0}$ in $\cF^{uu}(q_{X,p_{0}})$. In particular, for any $q \in B_E(q_{X,p_0},\rho_3) \cap S_0$, we know that $\cF^{uu}(q_{X,p_{0}}) \cap \cF^{s}(q) \not=\emptyset.$
Let 
\begin{displaymath}
U = \bigcup_{|t| < \rho_3} \{ g_t q \st q \in B_E(q_{X,p_0},\rho_3/2)
\cap S_0 \}.
\end{displaymath}
\noindent
{\it Claim:} there exists $\rho_0> 0$ depending only on $\cKM$, and $\widetilde{L}_{0}$ depending on $\epsilon$ and $\cKM$ such that 
for $L>\widetilde{L}_{0}$,
\begin{equation}
\label{eq:onto:ball}
\pi(g_L U) \supset B_{\cT}(\pi(\gamma(L)),\rho_0).
\end{equation}

This is straightforward (in view of the non-uniform hyperboliclity as in
(\ref{one}) and (\ref{two})) but somewhat tedious argument. Let
\begin{displaymath}
V = \bigcup_{|t| < \rho_3} \{ g_{t+t(q)} f(q) \st q \in B_E(q_{X,p_0},\rho_3/2)
\cap S_0 \} \subset \cF^{u}(q_{X,p_0}).
\end{displaymath}
Then $V$ is relatively open as a subset of $\cF^{u}(q_{X,p_{0}})$. Let
$\partial 
V$ denote the boundary of $V$ viewed as a subset of $\cF^{u}(q_{X,p_{0}})$. 
Therefore we can choose $\rho_2 > 0$ depending only on $\cKM$ such
that for all $q' \in \partial V$, 
\begin{displaymath}
d_E(q',q_{X,p_0}) \ge \rho_2.
\end{displaymath}
By $(\ref{one})$ and the fact that $\partial V \subset
\cF^{u}(q_{X,p_{0}})$, this implies that for some constant $\rho_2'$,
(depending only on $\cKM$), all $q' \in \partial V$ 
and all $L > 0$, 
\begin{equation}
\label{eq:gLqprime}
d_E(g_L q', \gamma(L)) \ge \rho_2'.
\end{equation}

Note that $U \subset \bigcup_{t \in \reals} g_t S_0 \isom \reals
\cross S_0$. Let $\partial U$ denote the boundary of $U$ viewed as a
subset of $\reals \cross S_0$. Suppose $q_1 \in \partial U$.
We may write $q_1 = g_t q$ for some $q \in S_0$. Then the fact that
$q_1 \in \partial U$ implies that either $d_E(q,q_{X,p_0}) = \rho_3/2$
or $|t| = \rho_3$. In either case, let $q_2 = g_{t+t(q)} f(q)$. Then
$q_2 \in \partial V$, $q_1$ and $q_2$ are on the same leaf of
$\cF^{ss}$, and
\begin{displaymath}
d_E(q_1,q_2) \le C,
\end{displaymath}
where $C$ depends only on $\cKM$. 
Hence, by $(\ref{one})$, we have 
 \begin{equation}\label{basic}
 d_E(g_L q_1, g_L q_2) \le C \cdot C_{0}.
 \end{equation}
 
In order to prove the claim, we show that there exists $L_0$ (depending only on $\cKM$ and $\epsilon$) such
that for $L > L_0$,
\begin{equation}
\label{claim}
d_E(g_L q_1, \gamma(L)) \ge \rho_2'/2.
\end{equation}
Suppose that $(\ref{claim})$ fails. Then:
\begin{itemize}
\item by $(\ref{basic})$, $d_{E} (g_{L}q_{2},\gamma(L)) \leq C_{2},$ where $C_2$ only depends on $\cKM.$ On the other hand, we can choose 
$|t_{0}| \leq \rho_{3}$ such that $g_{t_{0}}q_{2} \in \cF^{ss}(q_{X,p_0}).$
 Using $(\ref{two})$ for $g_{L+t_{0}}q_{2} \in \cF^{ss}(\gamma(L)) $, we get that 
$$|\{ t | \; t \in [0,L] , \cl(g_{t}q_{2}) > \epsilon' \}| > \epsilon_{0} s, $$
where $\epsilon', \epsilon_{0}$ only depend on $\cKM,$ and $\epsilon.$
\item Now we can apply $(\ref{two})$ for $q_{1} \in \cF^{ss}(q_{2})$, and we get $d_E(g_L q_1, g_L q_{2}) \leq C_{0} e^{-\na L} d_{H}(q_{1},q_{2}).$ Therefore there exists $\widetilde{L}_{0}$ such that if $L >\widetilde{L}_{0}$, $d_E(g_L q_1, g_L q_{2}) \leq
  \rho_2'/2$. Then, in view of (\ref{eq:gLqprime}), $d_{E}(g_{L}
  q_{2},\gamma(L)) \geq \rho_2'$ . So for all $q_1 \in \partial U$ and all $L > \widetilde{L}_0$, $ d_E(q_L q_1, \gamma(L)) \ge \rho_2'/2.$
\end{itemize}
By ${\bf P2}$, there exists $\rho_0$ depending only on $\cKM$ such
that for $L> \widetilde{L}_{0} $ and all $q_1 \in \partial U$, 
\begin{equation}
\label{eq:dct:lower:bound}
d_{\cT}(\pi(g_L q_1), \pi(\gamma(L))) \ge \rho_0.
\end{equation}
Let $\phi: U \to \cT_g$ denote the map $\phi(q) = \pi(g_L(q))$. Then
$\phi$ is continuous, and $\phi(q_{X,p_0}) = \pi(\gamma(L))$. Now in
view of $(\ref{eq:dct:lower:bound})$, $(\ref{eq:onto:ball})$ holds.
In other words, if $d_{\cT}(p_{0},p_{1}) < \rho_0,$ and $L>
\widetilde{L}_0$, $\cF^{uu}(q_{X,p_{0}})\cap \cF^{s}(q_{X,p_{1}}) =q$
is well defined, and $d_{E}(q, q_{X,p_{0}}) < \rho_{3}$. On the other
hand, $d_{H}(g_{L}q_{X,p_{0}}, g_{L}q) <C,$ where $C$ only depends on
$\cKM$. Therefore, we can use ${\bf P4}$ to get
$$d_{E}(q, q_{X,p_{0}}) < C_{0} e^{-\na L} C,$$
where $C_{0},\na>0$ depend on $\epsilon, \nb$ and $\cKM$. This shows that we can find $L_{0}>\widetilde{L}_0$ (depending on $\epsilon, \rho,$ and 
$\cKM$) such that if $L>L_{0}$ we have 
$$d_{E}(q, q_{X,p_{0}}) <\rho.$$
\qed\medskip\\
\noindent
{\bf Proof of Lemma $\ref{lemma:hyperbolic:nbhd}$.}
It is enough to show that given $\epsilon>0$ and $1 > m > m' \geq 1/2$, there exist $c_{0}>0$
(depending only on $\cKM$, $m$ and $m'$) and $L_{0}>0$ (depending only on
$\epsilon$, $m, m'$ and $\cKM$) such that if $x,p_{1} \in \cKT$,
$d_{\cT}( p_{1},p_{2}) < c_{0},$ for $i=1,2,$ $d_{\cT}(X,p_{i})=r_{i}> L_{0}$, and
$$\{ s\; | 0< s \leq r_{1}\; ,\; \cl(g_{s}(q_{X,p_{1}})) \geq \epsilon \} > r_{1} m ,$$
then 
$$\{ s\; | 0< s \leq r_{2}\; , \; \cl(g_{s}(q_{X,p_{2}}))\geq \epsilon/2 \} > r_{2} m'.$$
%Let $\gamma_{i}=\RR(q_{X,Y_{i}}),$ $\eta_{i}=\II(q_{X,Y_{i}}).$ 
%Also, let $\mathcal{K}=\mathcal{K}_{\epsilon},$ and $\mathcal{K} \subset \mathcal{K}'=\mathcal{K}_{\epsilon/2}.$
 Let $q_{1}=q_{X,p_{1}},$ and (as in the previous lemma) let $q=f(q_{X,p_{2}})$ be the unique point in 
 $\cF^{uu}(q_{1}) \cap \cF^{s}(q_{X,p_{2}}).$ 
 Now let $q_{2}$ be a quadratic differential of area $1$ on the geodesic joining $X$ to $p_{2}$ such that $q_{2} \in \cF^{ss}(q)$; in particular, we have 
\begin{equation}\label{all}
q \in \cF^{uu} (q_{1}),\;\; q \in \cF^{ss}( q_{2}).
\end{equation} 
By ${\bf P1}$, we can choose $\epsilon_{0}>0$ such that if $d_{E}(q',q'')<\epsilon_{0},$ and $\cl(q')>\epsilon$ then $\cl(q'')>\epsilon/2.$ 
Then by Lemma $\ref{lemma:compactness}$, ${\bf P4}$, and ${\bf P2}$, we can choose $c_{0}<1/2$ (depending only on $\cKM$) and $L_{0}>0$ such that for $d_{\cT}(X,p_{1})=r >L_{0},$ and $d_{\cT}(p_{1},p_{2})< c_{0}$ the following hold:

\begin{enumerate}
\item  As in ${\bf P4}$ (for $\nb=m)$, there exist $C_{0}>0$ and $\na>0$ such that $(\ref{two})$ holds. Then $L_{0}$ should satisfy:
$$C_{0} e^{-\na L_{0}/6} < \epsilon_{0}/4,$$
\item Also, by Lemma $\ref{lemma:compactness}$, if $L_{0}$ is large enough $$d_{E}(q_{1},q) < \epsilon_{0}/4C_{0}\; , \;\;\; d_{E} (q, q_{2}) \leq \epsilon_{0}/4C_{0} ,$$
and
$$d_{E}(g_{r}q_{1},g_{r}q_{2}) \leq 1.$$
\end{enumerate}
As a result, from $(\ref{one})$ we get
$$d_{E}( g_{r}q, g_{r} q_{2}) < \epsilon_{0}/4\; \mbox{and}\; d_{E}(g_{r}q_{1},g_{r}q)<C_{1},$$
where $C_{1}$ only depends on $\cKM,$ and $\epsilon.$

Consider the map between the points on the geodesic $[xp_{1}]$ to the points on $[xp_{2}]$ as follows:
$$ [xp_{1}] \rightarrow [xp_{2}]$$
$$ g_{s}q_{1}\rightarrow g_{s}q_{2}.$$

We can choose, $0\leq s_{0} \leq r$ such that 
$|\{ t\; | s_{0}< t <r\; , \cl(g_{t}q_{1}) \geq \epsilon\} | =r (m-m'), $ and let    
$$\mathcal{A}=\{ t\; | 0<t < s_{0}\; ,\cl(g_{s}q_{1}) \geq \epsilon\}.$$
It is easy to check that $|\mathcal{A}|> r m'. $
We claim that for $s\in \mathcal{A},$ we have 
\begin{equation}
\label{c1}
 d_{E}(g_{s}q_{1}, g_{s}q_{2}) \leq \epsilon_{0}.
 \end{equation}
 This is because:
\begin{itemize}
\item For any $s \in \mathcal{A}$, $(\ref{con})$ holds for $g_{r}q$ and $g_{r}q_{1}$ and the interval $(0, r-s)$. Hence, by $(\ref{two})$, 
$$d_{E}(g_{s}q_{1}, g_{s}q)=d_{E} (g_{s-r}( g_{r}q_{1}), g_{s-r}(g_{r}q))< \epsilon_{0}/4; $$ 
%and also $g_{s}q \subset \mathcal{K};$
\item by $(\ref{one})$, $$d_{E}(g_{s}q, g_{s}q_{2})< \epsilon_{0}/4;$$
\item Finally, since $q_{1} \in \cF^{uu}(q)$, $d_{E}(g_{s}q_{1}, g_{s}q_{2}) \leq \epsilon_{0}.$
\end{itemize}
Hence, for $s \in \mathcal{A}$, $\cl(g_{s}q_{1})\geq \epsilon,$ $(\ref{c1})$, and ${\bf P1}$ imply that 
$\cl(g_{s}q_{2}) \geq \epsilon/2.$
\qed\medskip

\subsection{Proof of Theorem~\ref{theorem:staying:in:principal:stratum}}
Choose $\theta_1 > 0$ \mc{how}.
% WHY IS THIS 1 here ??????????????????????????????????????????????Let $\cKM_1 \subset \cM_g$ be such that
Let $\cKM \subset \cM_g$ be such that
Theorem~\ref{theorem:staying:in:think:part} holds for $\cKM $,
and $\theta_1=\theta$. Let $\cKQ_2$ be a compact subset of 
$\cQ\cT_{g}(1,\dots,1)/\Gamma_{g}$ such that $\cKQ_2 \subset \gp^{-1}(\cKM)$,
\mc{fix notation} and let $\cKQ_3$ be a 
subset of the interior of $\cKQ_2$. We may choose these sets so that
$\mu(\cKQ_3) > (1/2)$, where $\mu$ is defined in \ch \ref{qd}. We also choose $\cKQ_2$ and $\cKQ_3$ to be
symmetric, i.e if $q \in \cKQ_2$ then $-q \in \cKQ_2$ (and same for $\cKQ_3$).
Then there exists $\epsilon > 0$ such that for $X \in \cKQ_3$,
$\cl(X) > \epsilon$. Let $c_0$ be as in
Lemma~\ref{lemma:hyperbolic:nbhd}. 
We now choose a $(c_1,c_2)$ separated net $\cN$ on $\cT_g$, 
which $c_1 < c_0$, $c_2 < c_0$. We may assume that $\cN \cap
\gp^{-1}(\cKM)$ is invariant under the action of the mapping class 
group. 

Suppose $X \in \cT_g$, and as before, let $S(X)$ denote the set of area $1$ holomoprhic quadratic differentials on $X$.
Let 
\begin{displaymath}
U(X,T) = \{ q \in S(X) \st |\{ t \in [0,T] \st g_t q \in \cKQ_2^c \}| >
(1/2)T \}.
\end{displaymath}
Let 
\begin{displaymath}
V(X,T) = \bigcup_{0 \le t \le T} \pi(g_t U(X,t)),
\end{displaymath}
\mc{worry about t vs T at end}
so that $V(X,T)$ is the subset of
$B_{\cT}(X,T)$ consisting of points $Y \in B_{\cT}(X,T)$ such that the geodesic from $X$ to $Y$ spends more than half the time
outside $\cKQ_2$.

By \cite[Theorem 2.6]{ABEM}, 
for any $\theta_1 > 0$, there exists $T > 0$ such that for any $\tau >
T$ and any $X \subset \cN \cap \gp^{-1}(\cKM)$, 
\begin{displaymath}
\mm(\Nbhd_{c_2}(V(X, \tau) \cap \gp^{-1}(\cKM))) \le \theta_1 e^{h \tau},
\end{displaymath}
where $\Nbhd_a(A)$ denotes the set of points within Teichm\"uller
distance $a$ of the set $A$. Then, since $\cKM$ is compact and
$\theta_1$ is arbitrary, this
implies that for any $\theta_2 > 0$ there exists $T > 0$ such that for
any $\tau > T$ and any $X \in \cN \cap \gp^{-1}(\cKM)$, 
\begin{displaymath}
|\cN \cap V(X,\tau) \cap \gp^{-1}(\cKM) | \le \theta_2 e^{h \tau}.
\end{displaymath}
By the compactness of $\cKM$ and \cite[Theorem 1.2 and Theorem 5.1]{ABEM}
there exists $C_1 > 1$ such that for $\tau$ sufficiently
large and any $X \in \cN \cap \gp^{-1}(\cKM)$, 
\begin{displaymath}
C_1^{-1} e^{h\tau } \le |\cN \cap B_{\cT}(X,\tau)| \le C_1 e^{h\tau}.
\end{displaymath}
Thus, for any $\theta_3 > 0$ there exists $T > 0$ such that for $\tau > T$,
\begin{equation}
\label{eq:most:netpoints:are:in:good:sector}
|\cN \cap V(X,\tau) \cap \gp^{-1}(\cKM) | < \theta_3 | \cN \cap B_{\cT}(X,\tau)|.
\end{equation}
From now on we assume that $\tau$ is sufficiently large so that
(\ref{eq:most:netpoints:are:in:good:sector}) holds.

Let $\cKM_1' = \Nbhd_{c_2}(\cKM)$, and 
let $\cG(R)$ denote the set of closed geodesics in $\cM_g$ of length
at most $R$, and let $\cG^{\cKM_1'}(\theta_3,R) \subset \cG(R)$ 
denote the subset which contributes to $N^{\cKM_1'}(\theta_3,R)$.
In view of Theorem~\ref{theorem:staying:in:think:part}, it is enough
to show that there exists $\delta_0 > 0$ such that 
for $R$ sufficiently large, 
\begin{equation}
\label{eq:GsetminusGK1}
|\cG(R) \setminus \cG^{\cKM_1'}(\theta_3,R)| \le e^{(h-\delta_0) R}. 
\end{equation}
As in \S\ref{sec:thin:part}, 
we associate a random walk trajectory $\Phi(\gamma)$ to each closed
geodesic $\gamma \in \cG(R)$. \mc{worry about choosing lifts}
Let $\cP_1(R) = \Phi(\cG(R)\setminus \cG^{\cKM_1'}(\theta_3,R))$ 
denote the set of resulting trajectories. 
Note that by construction, every trajectory in $\cP_1$ 
spends at most $\theta_3$ fraction of the time outside $\cKM$. 

Suppose $\Lambda = (\lambda_1,\dots, \lambda_n)$ be a trajectory in
$\cP_1(R)$. Let $J(\Lambda)$ denote the number of $j$, $1\le j \le n$
such that 
\begin{equation}
\label{eq:bad:j}
\lambda_j \in \gp^{-1}(\cKM^c) \quad\text{ or }\quad \lambda_{j+1} \in
V(\lambda_j,\tau) \cup \gp^{-1}(\cKM^c).
\end{equation}
Let $\cP_2(R) \subset \cP_1(R)$ denote the trajectories for which 
$J(\Lambda) < (\theta/2) n$. Then, as long as $\theta_3$ is chosen
sufficiently small \mc{explain}, the Law of Large Numbers implies
that there exists $\delta_1 > 0$ such that for $n=R/\tau$ 
sufficiently large, 
\begin{displaymath}
|\cP_1(R)\setminus \cP_2(R)| \le e^{(h-\delta_1)R}.
\end{displaymath}
Since every trajectory in $\cP_1(R)$ intersects the compact set
$\cKM$, for any $\Lambda \in \cP_1(R)$, the cardinality
of $\Phi^{-1}(\Lambda)$ is bounded by a constant $C$ depending only on
$\cKM$. \mc{maybe give reference}
Thus,
\begin{displaymath}
|\Phi^{-1}(\cP_1(R)\setminus \cP_2(R))| \le C e^{(h-\delta_1)R}.
\end{displaymath}
To complete the proof we will show that there exists a compact set
$\cKQ \subset \cQ\cM_{g}(1,\dots,1)$ 
such that any geodesic $\gamma \in \Phi^{-1}(\cP_2(R))$ spends at
least $1-\theta$ fraction of the time in 
$\cKQ$. Indeed, suppose $\Lambda = (\lambda_1,
\dots, \lambda_n) \in \cP_2(R)$ and $\gamma \in \Phi^{-1}(\Lambda)$. Then there exist points $q_j =
\gamma(t_j)  \in \cQ\cM_g$ such that $q_j \in \gamma$ and for $1 \le j \le n$,
$d_{\cT}(\pi(q_j), \lambda_j) \le c_2 < c_0$. 
Now suppose $j \not\in J(\Lambda)$. 
Then by Lemma~\ref{lemma:hyperbolic:nbhd}, there exists $\epsilon' >
0$ depending only on $\cKM$ and $\epsilon$, and a point
$t_j<t'<t_{j+1}$ such that $\cl(\gamma(t')) >\epsilon'$. But then 
for all $t_j < t < t_{j+1}$, $\cl(\gamma(t)) \ge e^{-\tau}
\epsilon'$, so the entire geodesic $\gamma$ spends at least 
$(1-\theta)$-fraction of the time in the compact set
$\{ q \in \cQ\cM_g \st \cl(q) > e^{-\tau} \epsilon'\}$.
\qed\medskip

\section{A Closing Lemma}\label{closing}
In this section, we use the properties of the geodesic flow discussed in 
\S \ref{rec} to prove the following closing lemma:
\begin{lemma}[Closing Lemma]
\label{lemma:closing}
Let $\cKQ$ be a compact subset of $\cQ\cM_g(1,\ldots,1)$ consisting of non-orbifold points of $\cQ\cM_g$. 
Given $q \in \cQ\cM_g(1,\ldots,1)$ and $\epsilon > 0$, there exist constants $L_0 > 0$, $\epsilon' > 0$, 
and open neighborhoods $U \subset U' $ of $q$ with the following property. 

Suppose that $\gamma: [0,L] \to \cQ\cM_g$ is a geodesic segment (parametrized by
arclength) such that 
\begin{itemize}
\item [{\rm (a)}] $\gamma(0), \gamma(L) \in U$ and 
\item [{\rm (b)}] \begin{displaymath}
|\{ t \in [0,L] \st \gamma(t) \in \cKQ \}| > (1/2) L,
\end{displaymath}
\item [{\rm (c)}] $L > L_0.$ 
\end{itemize}
Let  $\gamma_{1}$ be the closed path in $\cM_{g}$ which is the union of $\gamma$ and the geodesic segment from $\pi(\gamma(L))$ to $\pi(\gamma(0))$ in $\pi(U).$
Then there exists a unique closed geodesic $\gamma'$ with the following
properties:
\begin{itemize}
\item[{\rm (a')}] $\gamma'$ and $\gamma_{1}$ have lifts in $\cT_{g}$ which stay $\epsilon$ close (with respect to the Teichm\"uller metric on $\cT_{g})$.
\item[{\rm (b')}] The length of $\gamma'$ is within $\epsilon$ of $L$, 
\item[{\rm (c')}] $\gamma'$ passes through $U'$.
\end{itemize}
\end{lemma}
\noindent
{\bf Remark.} A stronger version of this statement can be found in \cite{Hamenstadt:BM} ($\ch 6$).

\bold{Outline of Proof.} 
Consider the stable and unstable foliations for the geodesic flow as in \ch \ref{CNP}. Our goal is to show that if 
$U$ is small enough, the first return map on these foliations will define a contraction with respect to the modified Hodge distance
 function. As a result, we find a fixed point for the first return map in $U'$ (which is the same as a closed geodesic going through $U'$).  

From Lemma~\ref{lemma:hyperbolic:nbhd}, we can find open neighborhoods $U \subset U'' \subset U'$ of $q$, and $L_0>0$ 
such that the following properties hold :
\begin{itemize}
\item 
if $\gamma: [0,L] \to \cQ\cM_{g}$ satisfies properties $(a)$, $(b)$ and $(c)$ then 
in view of the hyperbolicity
statement (\ref{two}) the time $L$ geodesic flow restricted
to the neighborhood $U''$ expands along the leaves of $\cF^{uu}$
and contracts along the leaves of $\cF^{ss}$, in the metric $d_H$, 
\item for any $q_{1}, q_{2} \in U',$ if $q_1 \in \cF^{ss}(q_2)$ or $q_1 \in \cF^{uu}(q_2)$ then
$d_{H}(q_{1},q_{2}) \leq \epsilon.$
\end{itemize}
We can apply the contraction mapping principle to $\cF^{ss}$ to find 
$q_{0} \in U'$ such that $g_{L}(q_{0})\in \cF^{uu} q_{0}.$
Now we can consider the first return map of the map $g_{-t}$ on $\cF^{uu}(q_0)$. Again, we can use 
 (\ref{two}), Lemma~\ref{lemma:hyperbolic:nbhd} and the contraction mapping principle to find
a fixed point for the geodesic flow in $U'$ .
\qed\medskip

\subsection{Remark on proof of Theorem~\ref{theorem:asymp:allgeodesics}}

Note that by the bound proved in Theorem \ref{theorem:staying:in:principal:stratum}, we only need to consider the 
set of closed geodesics going through a fixed compact subset of $\cM_{g}.$
We remark that in Lemma $\ref{lemma:closing}$ if we remove the assumption that $\cKQ$ consists of non-orbifold points then there are at most $c_{0}$ closed geodesics satisfying conditions ${\rm (a')},$ ${\rm (b')}$ and ${\rm(c')},$ where $c_{0}$ is a constant depending only on $g$.

Roughly speaking, since 
\begin{itemize}
\item by Theorem \ref{VM}, the geodesic flow on $\mathcal{Q}^{1}\mathcal{M}_{g}(1,\ldots,1)$ is mixing, and  
\item on a fixed compact subset of $\mathcal{Q}^{1}\mathcal{M}_{g}(1,\ldots,1)$ the geodesic flow is uniformly hyperbolic, 
\end{itemize} 
the derivation of Theorem \ref{theorem:asymp:allgeodesics} from Lemma \ref{lemma:closing} can be done 
following the work of Margulis \cite{Margulis:thesis}. See also $\S 20.6$ in \cite{Katok:Hasselblat}. We skip the argument since it has been done
carefully in $\ch 5$ and $\ch 6$ of  \cite{Hamenstadt:BM}.

% For details in the Teichm\"uller space setting see

\section{Appendix: Proof of Proposition~\ref{prop:measures:of:balls}}

Given $X$ and $Y$ in $\cT_g$ and $\tau > 0$, let $F_\tau(X,Y)$
denote the intersection of $B_{\cT}(X,\tau)$ with
$\Gamma_{g}\cdot Y$. Proposition~\ref{prop:measures:of:balls} will follow
from the following:
\begin{prop}
\label{prop:toprove}
Given $\epsilon>0,$ there exists $\tau_{0}>0$ such that for any $\tau>\tau_{0}$ and $X,Y \in \te_{g}$
we have $$ | F_{\tau}(X,Y)| \leq e^{(h+\epsilon) \tau} G(Y)^{2},$$
where $G$ is defined in $(\ref{GG}).$
\end{prop}

\bold{Proof of Proposition~\ref{prop:measures:of:balls} assuming
Proposition~\ref{prop:toprove}.}
Suppose $X \in \cT_g$ is arbitrary, and suppose 
$Y \in \widetilde{\cN}$. Then, by Proposition~\ref{prop:toprove},
for $\tau$ sufficiently large,  
\begin{displaymath}
|B_{\cT}(X,\tau) \cap \Gamma_{g} Y| \leq e^{(h+\epsilon/2) \tau} G(Y)^{2}.
\end{displaymath}
Note that by Theorem $\ref{theorem:minsky:ext}$ (see the calculation in $(\ref{simplebound})$) 
\begin{equation}
\label{eq:spread}
| \Gamma_{g} \cdot Y \;\cap B_{\cT}(Y,c_{2})| \geq c\; G(Y)^{2}. 
\end{equation}
Thus, since distinct points in $\cN$ are are least $c_{2}$ apart, 
\begin{displaymath}
|B_{\cT}(X,\tau) \cap \Gamma_{g} Y \cap \cN| 
\leq c \; e^{(h+\epsilon/2) \tau}.
\end{displaymath}
In view of (\ref{eq:netpoints}), this implies that
\begin{displaymath}
|B_{\cT}(X,\tau) \cap \cN| 
\leq c' \tau^{6g-6} e^{(h+\epsilon/2) \tau} \le e^{(h+\epsilon)\tau}. 
\end{displaymath}
Now this implies Proposition~\ref{prop:measures:of:balls} in view of
Lemma~\ref{lemma:uniform:vol:ball}. 
\qed\medskip

Given $X, Y \in
\cT_{g}$, and $\mathcal{B} \subset \cC_{X}$ let
$$F_{R}(X,Y,\mathcal{B})=  \{\gg\cdot Y\; | \gg \in \Gamma_{g}, d_{\cT}(X, \gg \cdot Y) \leq R, \; \mathcal{B} \subset \cC_{X}\cap \cC_{\gg \cdot Y}\} \subset F_{R}(X,Y).$$
Here, as in \ch \ref{sec:back}, given $Z \in \cT_{g}$ the set $\cC_{Z}$ consists of simple closed curves of extremal length $\le \epsilon_0^2$ on $Z$.

In general, we have:
\begin{theorem}
\label{theorem:uniform:upper:bound}
Given $X ,Y \in \cT_g$ 
\begin{equation}\label{weknow}
|F_{R}(X,Y,\mathcal{B})| \leq  C_{1} R^{3g-3} e^{(h-2|\mathcal{B}|) R} G(X) G(Y), 
\end{equation}
 where $C_{1}$ only depends on $g.$
\end{theorem}
In the proof we use the following lemma from \cite{ABEM}: 
\begin{lemma}\label{Mmain}
Let $\alpha=\{\alpha_{1},\ldots,\alpha_{3g-3}\}$ be a bounded pants
decomposition on $X$ and suppose $Y_{0} \in B_{\cT}(X,R)$.
If $h_{\alpha_{1}}^{m_{1}} \cdots h_{\alpha_{3g-3}}^{m_{3g-3}} (Y_{0}) \in B_{\cT}(X,R)$ then for $1 \leq i\leq 3g-3$
$$|m_{i}| \cdot \sqrt{\e_{\alpha_{i}}(Y_{0})} \leq C \frac{e^R}{\sqrt{\e_{\alpha_{i}}(X)}}.$$
Here $C$ is a constant which only depends on $g.$
\end{lemma}
\noindent
{\bf Proof of Theorem \ref{theorem:uniform:upper:bound}.} 
To simplify the notation, given $X, Y \in \cT_{g}$ define $$D_{\alpha}(X,Y)=\log \left (\sqrt{\frac{\e_{\alpha}(X)}{\e_{\alpha}(Y)}}\right).$$
Choose bounded pants decompositions $\mathcal{P}_{X}=\{\alpha_{1},\ldots, \alpha_{3g-3}\}$ and $\mathcal{P}_{Y}=\{\beta_{1},\ldots, \beta_{3g-3}\}$ on $X$ and $Y$. 
Without loss of generality, we can assume that 
$$ \alpha_{1}=\beta_{1},\ldots, \alpha_{b}=\beta_{b} \in \mathcal{B},$$
 where $b=|\mathcal{B}|.$

Given ${\bf r}=(r_{1},\ldots, r_{3g-3}), {\bf s}=(s_{1},\ldots,s_{3g-3}) \in {\Bbb Z}^{3g-3}$ define
$$M_{R, ({\bf r},{\bf s})}(X,Y,\mathcal{B})= M_{R}(X,Y,\mathcal{B}) \cap \{ Z= \gg \cdot Y |\; D_{\alpha_{i}}(X,\gg Y) \in [s_{i},s_{i}+1] , D_{\beta_{i}}(\gg^{-1} X,Y) \in [r_{i},r_{i}+1] \}.$$
By the definition, if $M_{R, ({\bf r,\bf s})}(X,Y,\mathcal{B}) \not =\emptyset $ then for $1 \leq i \leq 3g-3$ we have $|r_{i}|, |s_{i}| \leq R. $\\

\noindent
{\it Claim.}
Given ${\bf r}, {\bf s} \in  {\Bbb Z}^{3g-3}$ with $|r_{i}|, |s_{i}| \leq R$ 
\begin{equation}
\label{c}
|M_{R, ({\bf r,\bf s})}(X,Y,\mathcal{B}) | \leq  e^{(h-2|\mathcal{B}|) R} G(X) G(Y).
\end{equation}
In order to prove this claim, fix $Y_{0}= \gg_{0} \cdot Y \in M_{R, ({\bf r,\bf s})}(X,Y,\mathcal{B}),$
and consider the pants decomposition $\alpha= \cup_{i=1}^{3g-3} \alpha_{i}.$
Note that given $Z= \gg\cdot Y \in M_{R, ({\bf r,\bf s})}(X,Y,\mathcal{B})$,  $$\alpha(Z)=\gg_{0} \gg^{-1} (\mathcal{P}_{X})=\cup_{i=1}^{3g-3} \gg_{0} \gg^{-1} (\alpha_{i}) $$ defines a pants decomposition. Moreover, we have:

\begin{itemize}
\item $\gg_{0} \mathcal{P}_{Y}$ is a bounded pants decomposition on $Y_{0}=\gg_{0} \cdot Y.$
\item Also, 
\begin{equation}\label{BB}
\sqrt{\e_{Y_{0}}(\alpha(Z))}= \sqrt{\e_{Z}(\alpha)}= O(e^{R}).
\end{equation}
\item 
Note that $$\cup_{i=1}^{b} \alpha_{i} \subset \alpha(Z),$$
and hence for $i \leq b$
\begin{equation}\label{Y}
i( \gg_{0} \beta_{i}, \alpha(Z))=0.
\end{equation}
\item 
On the other hand, since $Z= \gg \cdot Y$ is in $M_{R, ({\bf r,\bf s})}(X,Y,\mathcal{B})$
we have:
$$   \frac{\sqrt{\e_{\beta_{j}}(\gg^{-1} X)}}{  \sqrt{\e_{\beta_{j}}(Y)}} =O(e^{r_{i}}).$$
As a result, we get 
$$i(\gg_{0} \gg^{-1}(\alpha_{i}) , \gg_{0} \beta_{j})=i(\alpha_{i}, \gg (\beta_{j})) \leq \sqrt{\e_{\alpha_{i}}(X)} \; \sqrt{\e_{\beta_{j}}(\gg^{-1}X)}=O(e^{r_{j}}),$$
and hence
\begin{equation}\label{YY}
 i(\alpha(Z), \gg_{0}\cdot \beta_{j})=O(e^{r_{j}}). 
 \end{equation}
\item By $(\ref{BB})$ Theorem $\ref{theorem:minsky:ext} $ gives rise to an upper bound on twisting parameter of $\alpha(z)$ around $\gg_{0} \beta_{j}$ on $Y_{0}$:
\begin{equation}\label{D}
 tw_{Y_{0}}(\alpha(Z), \gg_{0} \beta_{j}) = tw_{Y}(\gg^{-1}(\alpha), \beta_{j}) \leq \frac{e^{R}}{\sqrt{\e_{\beta_{j}}(Y)}}.
 \end{equation}
 
Now, in view of the Dehn-Thurston parametrization of multicurves, $(\ref{Y})$, $(\ref{YY})$ and $(\ref{D})$ imply:
\begin{equation}\label{y}
| \{\alpha(Z)\;| Z \in M_{R, ({\bf r,\bf s})}(X,Y,\mathcal{B})\}| \leq e^{r_{b+1}+\ldots r_{n}} e^{(\frac{h}{2}-b)R} 
\prod_{i=1}^{3g-3}\frac{1}{\sqrt{\e_{\beta_{i}}(Y)}} 
\end{equation}

\item On the other hand given $Z_{1}=\gg_{1}\cdot Y, Z_{2}=\gg_{2}\cdot Y \in \Gamma_{g}\cdot
y$, we have $\alpha(Z_{1})=\alpha(Z_{2})$ if and only if $\gg_{1}^{-1}\; \gg_{2} (\mathcal{P}_{X})=\mathcal{P}_{X}.$
If  $\gg_{1}^{-1}\; \gg_{2}(\alpha_{i})=\alpha_{i}$ for $1\leq i \leq 3g-3$ then 
there are
$m_{1},\ldots m_{3g-3} \in{\mathbb Z}$ such that
$$\gg_{1}=h_{\alpha_{1}}^{m_{1}} \cdots h_{\alpha_{3g-3}}^{m_{3g-3}}
\cdot \gg_{2}.$$  

Since for any $Z \in M_{R, ({\bf r,\bf s})}(X,Y,\mathcal{B})\ $
$$ \frac{1}{\sqrt{\e_{\alpha_{i}}(Z)}}=O( \frac{1}{e^{s_{i}}}),$$ 
from Lemma \ref{Mmain} we get
\begin{equation}\label{d}
| \{ Z \in M_{R, ({\bf r,\bf s})}(X,Y,\mathcal{B})\; | \; \alpha(Z)=\alpha(Z_{0}) \}| =O \left(\frac{e^{(\frac{h}{2}-b)R}}{e^{|{\bf s}|}}\prod_{i=1}^{3g-3} \frac{1}{\sqrt{\e_{\alpha_{i}}(X) }}\right). 
\end{equation}
\end{itemize}

Therefore from $(\ref{y})$ and $(\ref{d})$ we get 
$$ |M_{R, ({\bf r,\bf s})}(X,Y,\mathcal{B}) | =O\left( e^{|{\bf r}|-|{\bf s}|} G(X)G(Y) e^{(h-2|\mathcal{B}|)R}\right) . $$ 
On the other hand, since $M_{R, ({\bf r,\bf s})}(X,Y,\mathcal{B})= M_{R, ({\bf s,\bf r})}(Y, X,\mathcal{B}),$ we have
$$  |M_{R, ({\bf r,\bf s})}(X,Y,\mathcal{B})|= O\left( e^{|{\bf s}|-|{\bf r}|} G(X) G(Y) e^{(h-2|\mathcal{B}|)R}\right).$$
Therefore, we get $(\ref{c})$.

Finally, we have
$$ M_{R}(X,Y,\mathcal{B})= \bigcup  M_{R, ({\bf r,\bf s})}(X,Y,\mathcal{B}), $$
where ${\bf r}=(r_{1},\ldots,r_{3g-3})$ and ${\bf s}=(s_{1},\ldots,s_{3g-3})$ with $|r_{i}|. |s_{i}| \leq R.$ Hence, in view of $(\ref{c})$ we get $(\ref{weknow})$.

\hfill $\Box$

\bold{Proof of Proposition~\ref{prop:toprove}.}
To simplify the notation, given $X \in \cT_{g}$ let $${\bf e}_\alpha(X)=\sqrt{\e_{\alpha}(X)}.$$
As before, we say $\alpha$ is short on $X$ if
${\bf e}_\alpha(X) \leq \epsilon_{0},$ and let $\mathcal{C}_{X}$ denote the
set of all short simple closed curves on $X.$ 
Choose a $(c,2c)$ net $\cN$ in $\te_{g}$ as in $\ch \ref{netpoints}.$
%Given $X, Y \in
%\te_{g}$, and $\mathcal{B} \subset \mathcal{A}_{X}$ let
%$$F_{\tau}(X,Y,\mathcal{B})=  \{gY\; | g \in \Gamma, d_{\cT}(X,gY) \leq \tau, \mathcal{A}_{X}\cap \mathcal{A}_{gY}=\mathcal{B}\} \subset F_{\tau}(X,Y).$$
%Using the method in \cite[\S{10}]{ABEM}, we have
%\begin{equation}\label{weknow}
%|F_{\tau}(X,Y,\mathcal{B})| \leq  C_{1} e^{(h-|\mathcal{B}|) \tau} G(X) G(Y).
%\end{equation}

Then in view of (\ref{eq:spread}), given $Z \in \te_{g}$,
$\mathcal{B} \subset \mathcal{C}_{Z}$, $r>0$ and $W \in \widetilde{\cN}$,
the inequality $(\ref{weknow})$ implies that:
\begin{align}
\label{weknow2}
| B_{\cT}(Z,\tau) \cap \gp^{-1}(W) \cap \cN_{\mathcal{B}}| & \leq C_{1}
e^{(h-|\mathcal{B}|) r} r^{3g-3} G(Z) G(W)/ G(W)^{2} \notag \\
& \leq C r^{3g-3} e^{(h-|\mathcal{B}|)r} G(Z)/ G(W) ,
\end{align}
where $\cN_{\mathcal{B}}=\{Z_{1} \in \cN | \forall \alpha \in \mathcal{B}, {\bf e}_\alpha (Z_{1}) \leq \epsilon_{0}\} \subset \cN,$
and $C$ is a universal constant independent of $X$, $Z$ and $r$. 

We order the elements in $\mathcal{C}_{X}$, $(\alpha_{1},\ldots,\alpha_{k})$, so that  
$${\bf e}_{\alpha_{1}}(X)\geq \ldots \geq {\bf e}_{\alpha_{s}}(X) \geq e^{-\tau}> \ldots {\bf e}_{\alpha_{k}}(X).$$
Then $3g-3\geq k \geq s.$ Let $\tau_{0}=0$, and for $i\geq 1$, $$\tau_{i}=-\log({\bf e}_{\alpha_{i}}(X)).$$ Also for $1\leq i \leq s$, $d_{i}=\tau_{i}-\tau_{i-1},$
$m_{i}=h-k+i-1,$ $d_{s+1}=\tau-\tau_{s},$ and $m_{s+1}=h-k.$
Then it is easy to verify that 
\begin{enumerate}
\item 
$$\sum_{i=1}^{s+1} d_{i}= \tau, \;\; \sum_{i=1}^{s+1} m_{i} d_{i}=(h-k+s) \tau-\tau_{1}\ldots-\tau_{s}, G(X)=e^{\tau_{1}+\ldots+\tau_{k}};$$
\item By Theorem $\ref{ker}$ , for $1\leq i\leq s$, $\alpha_{i},\ldots,\alpha_{k}$ are all short in $B_{\cT}(X,\tau_{i})$. 
\item Since $\alpha_{s+1},\ldots \alpha_{k}$ stay short in $B_{\cT}(X,\tau)$, again by Kerckhoff's theorem, for any $Y \in B_{\cT}(X,\tau)$
\begin{equation}\label{vs}
e^{\tau_{s+1}+\ldots+\tau_{k}} \leq G(Y) e^{(k-s) \tau}.
\end{equation}

\end{enumerate}
Consider 
$$\mathcal{P}=\{(X=Z_{0}, Z_{1},\ldots,Z_{s}, \gg\cdot Y=Z_{s+1})\; | Z_{1},\ldots,Z_{s} \in \cN, \forall \; i , d_{\cT}(Z_{i},Z_{i-1})\leq d_{i}, \gg \in \Gamma_{g}\},$$
$$\mathcal{Z}=\{ (W_{1},\ldots, W_{s})\;| W_{i} \in \widetilde{\cN} \cap 
\gp(B_{\cT}(X,\tau))\},$$
and finally  
 $$ \mathcal{P}(W_{1},\ldots W_{s})=\{(X, Z_{1},\ldots,Z_{s}, \gg\cdot Y) \in \mathcal{P}\;| \pi(Z_{i})=W_{i}\}.$$ 
Note that for any $(X,Z_{1},\ldots,Z_{s}, Z_{s+1}=\gg \cdot Y) \in \mathcal{P}$, $Z_{i} \in B_{\cT}(X,\tau_{i}).$ 

On the other hand, since we can approximate a geodesic by points in the net $\cN$, we have
$$|F_{\tau}(X,Y)| \leq |\mathcal{P}|,$$
also by the definition, $$\mathcal{P} \subset \bigcup_{W\in \mathcal{Z}} \mathcal{P}(W).$$
By equation $(\ref{eq:netpoints})$
we have 
$$|\mathcal{Z}|=O(\tau^{c_{g}}) ,$$
where $c_{g} \leq (12g-12)s < (12g-12)^{2}.$

Let $Z_{0}=X.$ Then for a given $Z=Z_{i-1} \in \cN \cap
B_{\cT}(X,\tau_{i-1})$, $B_{\cT}(Z,d_{i}) \subset B_{\cT}(X,\tau_{i}).$ This implies
that $h-m_{i}$ simple closed curves (i.e, the curves in $ \mathcal{B}_{i}=\{
\alpha_{i},\ldots, \alpha_{k}$\}) are short on $B_{\cT}(Z,d_{i}).$

Now, since $d_{i} \leq \tau,$ equation $(\ref{weknow2})$ for $\mathcal{B}=\mathcal{B}_{i},$
$W=W_{i}$, and $Z=W_{i-1}$ implies that
\begin{equation}
\label{more} 
|\gp^{-1}(W_{i}) \cap \cN \cap B_{\cT}(Z_{i-1},d_{i})|= 
|\gp^{-1}(W_{i}) \cap \cN_{\mathcal{B}_{i}} \cap B_{\cT}(Z_{i-1},d_{i})|
 \leq C \tau^{3g-3} e^{ m_{i} d_{i}} G(W_{i-1})/G(W_{i}).
 \end{equation}

Given $(Z_{1},\ldots,Z_{s}) \in \mathcal{Z}$, we can apply $(\ref{more})$ for $(X,W_{1}),\ldots, (W_{s-1}, W_{s})$, and apply $(\ref{weknow})$ for $(W_{s}, Y)$. We get 

$$|\mathcal{P}(W_{1},\ldots,W_{s})| \leq C \tau^{3g-3} e^{m_{1} d_{1}} G(X)/G(W_{1}) \ldots C \tau^{3g-3} e^{m_{s} d_{s}} G(W_{s-1})/G(W_{s})\times $$ $$\times C \tau^{3g-3} e^{m_{s+1} d_{s+1}} G(W_{s}) G(Y) 
=C^{s} \tau^{(3g-3)s} \frac{e^{(h-k+s)\tau}}{e^{\tau_{1}+\ldots \tau_{s}}} G(X) G(Y) = $$
$$=C^{s} \tau^{(3g-3)s} e^{(h-k+s)\tau} \times e^{\tau_{s+1}+\ldots+ \tau_{k}} G(Y)
\leq C^{3g-3} \tau^{(3g-3)^{2}} e^{h \tau} G(Y)^{2}.$$
We are using $(\ref{vs})$ to obtain the last inequality.
Now we have, 
$$|\mathcal{P}| \leq |\mathcal{Z}|\; e^{h\tau} G(Y)^{2} \leq C' \tau^{c'_{g}} e^{h \tau} G(Y)^{2},$$
where $c'_{g}=c_{g}+(3g-3)^{2}= O(g^{2}),$ and $C'=C^{3g-3}.$
\hfill $\Box$

\end{document}